\documentclass[leqno,12pt,a4paper]{amsart}%
\usepackage{amsmath,amssymb}
\usepackage{amsfonts}%
\usepackage{color}

\newtheorem{Theorem}{\hspace{\parindent}\bf Theorem}[section]
\newtheorem{Lemma}[Theorem]{\hspace{\parindent}\bf Lemma}
\newtheorem{Proposition}[Theorem]{\hspace{\parindent}\bf Proposition}
\newtheorem{Corollary}[Theorem]{\hspace{\parindent}\bf Corollary}

\newtheorem{Remark}[Theorem]{\hspace{\parindent}\bf Remark}

\begin{document}

\title[Logarithmic corrections in Fisher-KPP type PMEs]{\textbf{Logarithmic corrections in Fisher-KPP type Porous Medium Equations}$^*$}
\author[Y. Du, F. Quir\'{o}s and M. Zhou]{Yihong Du$^\dag$, Fernando Quir\'{o}s$^\S$ and Maolin Zhou$^\ddag$}
\thanks{\hspace{-.5cm}$^*$This research was supported by the Australian Research Council, and by 
	projects MTM2014-53037-P  and MTM2017-87596-P (Spain).}
\thanks{\hspace{-.6cm}
$^\dag$School of Science and Technology, University of New England, Armidale, NSW 2351, Australia\\ 
{\sf  (e-mail: ydu@une.edu.au)}
\\
 $^\S$Departamento de Matem\'{a}ticas, Universidad
Aut\'{o}noma de Madrid, 28049 Madrid, Spain \\ 
{\sf (e-mail: fernando.quiros@uam.es)}
\\
$^\ddag$School of Science and Technology, University of New England, Armidale, NSW 2351, Australia\\ 
{\sf(e-mail: mzhou6@une.edu.au)}
}

\date{\today}

\begin{abstract}
We consider the large time behaviour of solutions to the porous medium equation with a Fisher-KPP type reaction term and  nonnegative,  compactly supported initial function in $L^\infty(\mathbb{R}^N)\setminus\{0\}$:
\begin{equation}
\label{eq:abstract}
\tag{$\star$}u_t=\Delta u^m+u-u^2\quad\text{in }Q:=\mathbb{R}^N\times\mathbb{R}_+,\qquad u(\cdot,0)=u_0\quad\text{in }\mathbb{R}^N,
\end{equation}
with $m>1$. It is well known that the spatial support of the solution $u(\cdot, t)$ to this problem remains bounded for all time $t>0$.
In spatial dimension one it is known that there is a minimal speed $c_*>0$ for which the equation admits a wavefront solution $\Phi_{c_*}$ with a finite front, and the associated traveling wave solution is asymptotically stable in the sense that if the initial function $u_0\in L^\infty(\mathbb R)$ satisfies $\liminf_{x\to -\infty} u_0(x)>0$ and $u_0(x)=0$ for all large $x$, then
$\lim_{t\to\infty} \left\{\sup_{x\in \mathbb R}|u(x,t)-\Phi_{c_*}(x-c_*t-x_0)|\right\}=0$ for some $x_0\in \mathbb R$.  In dimension one we can obtain an analogous stability result for the case of compactly supported initial data. In higher dimensions we show that $\Phi_{c_*}$ is still attractive, albeit that a logarithmic shifting occurs. More precisely, if the initial function in \eqref{eq:abstract} is additionally assumed to be radially symmetric, then  there exists a second constant $c^*>0$  independent of the dimension $N$ and the initial function $u_0$, such that
\[
\lim_{t\to\infty}\left\{\sup_{x\in\mathbb R^N}\big|u(x,t)-\Phi_{c_*}(|x|-c_*t+(N-1)c^*\log t-r_0)\big|\right\}=0
\]
for some $r_0\in\mathbb{R}$ (depending on $u_0$).
  If the initial function is not radially symmetric, then there exist $r_1, r_2\in \mathbb{R}$  such that  the boundary of the spatial support of the solution $u(\cdot, t)$ is contained in the spherical shell $\{x\in\mathbb R^N: r_1\leq |x|-c_* t+(N-1)c^* \log t\leq r_2\}$ for all $t\ge1$. Moreover, as $t\to\infty$, $u(x,t)$ converges  to $1$ uniformly in $\big\{|x|\leq c_*t-(N-1)c\log t\big\}$   for any $c>c^*$.
\end{abstract}

\maketitle

\vfill

\noindent{\makebox[1in]\hrulefill}\newline2010 \textit{Mathematics Subject
Classification.}
35K55, 
35B40. 
35K65, 
76S05. 

\newpage

\section{Introduction}
\setcounter{equation}{0}

The aim of this paper is to characterize the precise large time behaviour of solutions to the porous medium equation with a Fisher-KPP source term
\begin{equation}
\label{eq:main}
u_t=\Delta u^m+u-u^2\quad\text{in }Q:=\mathbb{R}^N\times\mathbb{R}_+,\qquad u(\cdot,0)=u_0\quad\text{in }\mathbb{R}^N,
\end{equation}
 where the initial function $u_0\not\equiv0$ is  bounded, nonnegative and compactly supported. We always assume that $m>1$. It is possible to replace $u^m$ and $u-u^2$ in \eqref{eq:main}  by more general functions of a similar type, but, for simplicity and clarity of presentation, we
only consider these special nonlinearities in this paper.

Problem \eqref{eq:main} arises from a variety of applications. For example, it was used by Gurtin and MacCamy~\cite{Gurtin-MacCamy-1977}  to describe the growth and spread of a spatially distributed biological population whose tendency to disperse depends on the population density. A similar model with $m=2$ (though with a Malthusian instead of a Fisher-KPP growth term),  had been introduced earlier  by Gurney and Nisbet~\cite{Gurney-Nisbet-1975}. See also~\cite{SanchezGarduno-Maini-1994} and the references therein for further background on this kind of models, and \cite{Perthame-Quiros-Vazquez-2014} for a recent application of problem~\eqref{eq:main} to the description of tumor growth.

Problem \eqref{eq:main} can be regarded as a generalization of the semilinear case $m=1$, which was introduced independently by Fisher~\cite{Fisher-1937}, to study some spreading questions in population genetics, and by Kolmogorov, Petrovsky and Piskunoff~\cite{Kolmogorov-Petrovsky-Piscounov-1937}, who considered also more general reaction nonlinearities, with a similar motivation.

The equation in~\eqref{eq:main} is degenerate when $u=0$, and it does not have in general classical solutions.  A function $u\in C(Q)$ is a \emph{weak solution} to problem~\eqref{eq:main} if it satisfies
\begin{equation}
\label{eq:weak.solution}
\int_{\mathbb{R}^N}u(t)\varphi(t)=\int_{\mathbb{R}^N}u_0\varphi(0)+\int_0^t\int_{\mathbb{R}^N}
\big[u\varphi_t+u^m\Delta\varphi+u(1-u)\varphi\big]
\end{equation}
for each $\varphi\in C^\infty(Q)$ with compact support in $x$ for every $t\ge0$. If the equality in~\eqref{eq:weak.solution} is replaced by \lq\lq$\le$'' (respectively, by \lq\lq$\ge$'') for $\varphi\ge 0$, we have a subsolution (respectively, a supersolution). Problem~\eqref{eq:main} has a unique bounded weak solution when the initial function is in $L^1(\mathbb{R}^N)\cap L^1(\mathbb{R}^N)$; see  Theorem 3 in \cite{Brezis-Crandall-1979}.  This solution satisfies moreover that $\nabla u^m \in L^2_{\rm loc}(Q)$;  see~\cite{Sacks-1983}, where $u_0$ is only required to be bounded. More general conditions on
$u_0$ can be found in \cite{dePablo-Vazquez-1991}, where a comparison principle between bounded subsolutions and supersolutions is also given.

A main feature of problem~\eqref{eq:main}, arising from its degeneracy at  $u=0$, is that its solution propagates with finite speed: 
if the initial function is compactly supported, the same is true for $u(\cdot,t)$ for all $t>0$. This can be seen  by a simple comparison 
with the solution of $u_t=\Delta u^m+u$ with the same initial function. Thus solutions to problem~\eqref{eq:main} have 
a \emph{free boundary}, 
separating the positivity region of $u$ from the region where $u$ vanishes. One of the main purposes 
of this paper is to describe precisely the large time behavior of the free boundary in the radial case, which then
provides  important information for the nonradial case via simple applications of the comparison principle.

The large time behavior of both the solution to~\eqref{eq:main} and its free boundary will be given in terms of a one-dimensional \emph{traveling wave}, that is, a solution to the equation in dimension $N=1$ of the form $\bar u_c(x,t)=\Phi_c(x-ct)$ for some \emph{speed} $c>0$ and \emph{profile} $\Phi_c$, which depends on $c$.

In order for $\bar u_c$ to be a weak solution of the equation, we look for profiles such that $\Phi_c\in C(\mathbb{R})$, $(\Phi_c^m)'\in L^2_{\rm loc}(\mathbb{R})$ and
$$
\int_{\mathbb{R}}\big((\Phi_c^m)'\varphi'+c \Phi_c\varphi'-\Phi_c(1-\Phi_c)\big)=0
$$
for all $\varphi\in C^\infty_{\rm c}(\mathbb{R})$.    We are interested in monotonic traveling wave profiles connecting the two equilibrium states 1 and 0:
\[
\Phi_c(-\infty)=1,\; \Phi_c(+\infty)=0,\; \Phi_c'\leq 0.
\]
 These are called, following~\cite{Fife-1978,Murray-1989}, \emph{wavefronts} (from 1 to 0) with speed $c$. Clearly
if $\Phi_c$ is a wavefront with speed $c$, then so is its translation $\Phi_c(\cdot+\bar\xi)$ for any fixed $\bar\xi\in\mathbb{R}$.

It is well known that for every $m > 1$ there exists a critical wave speed  $c_* = c_*(m) > 0$ such that there are
no wavefronts for \eqref{eq:main} with speed $c$ satisfying $0 < c < c_*$, while  there
exists exactly one wavefront (up to translation) for each speed $c\ge c_*$. If $c=c_*$, the corresponding wavefront is finite: there exists $\xi_*\in\mathbb{R}$ such that $0<\Phi_{c_*}(\xi)<1$ for $\xi<\xi_*$, while $\Phi_{c_*}(\xi)=0$ for $\xi\ge\xi_*$.
On the other hand, when $c>c_*$, the corresponding wavefronts are positive: $0<\Phi_c<1$ in $\mathbb{R}$. All these wavefronts satisfy $(\Phi_c^m)'\in C(\mathbb{R})$. See~\cite{Aronson-1980,Atkinson-Reuter-RidlerRowe-1981} for a proof of these facts using phase-plane analysis, and~\cite{Gilding-Kersner-2004} for an alternative proof using singular Volterra-type integral equations.

The wavefronts that describe the large time behaviour of solutions to \eqref{eq:main} with radial, compactly supported initial data are the ones with critical speed $c_*$. From now on we will denote by $\Phi_{c_*}$ the  wavefront of speed $c_*$ which is shifted so that $\xi_*=0$.  When $m=2$, both the critical speed and the corresponding wavefront are explicit, namely $c_*=1$ and $\Phi_{c*}(\xi)=\big(1-e^{\xi}\big)_+$; see~\cite{Aronson-1980, Newman-1980}. For the general case an explicit expression is no longer available, but we know (see \cite{dePablo-Vazquez-1991}) that
$\Phi_{c_*}\in C^\alpha(\mathbb{R})$, $\Phi_{c_*}\in C^\infty(\mathbb{R}_-)$, $\Phi_{c_*}'<0$ in $\mathbb{R}_-$, and
\begin{equation}
\label{eq:fb.condition.profile}
\lim_{\xi\to0^-} \frac{m}{m-1}\big(\Phi_{c_*}^{m-1}\big)'(\xi)=-c_*.
\end{equation}
This hints on the importance of the so called \emph{pressure} variable  $v=\frac{m}{m-1}u^{m-1}$. Indeed,~\eqref{eq:fb.condition.profile} implies that the velocity of the free boundary of the critical traveling wave $\bar u_{c_*}(x,t)=\Phi_{c_*}(x-c_*t)$ is given by the gradient of the pressure at this point. As we will see later, this is also true for general radial solutions. Notice that $\nabla u^m=u\nabla v$. Hence, from the point of view of mechanics, the velocity $V$ of the free boundary is given by Darcy's law, $V=-\nabla v$ at the free boundary, which is often a main ingredient in the modeling process.



\medskip

It turns out that $c_*$ can be regarded as the asymptotic spreading speed of the species. Indeed, if $u$ is a solution of~\eqref{eq:main} with a continuous and compactly supported nonnegative initial function $u_0\not\equiv 0$, then
\begin{equation}\label{spreading-rough}
\left\{\begin{array}{l}
\mbox{ $u(x,t)\to 1$ uniformly in $\{|x|<ct\}$ as $t\to\infty$, for any $c\in(0,c_*)$},\smallskip\\
\mbox{ $u(x,t)\to 0$  uniformly in $\{|x|>ct\}$ as $t\to\infty$, for any $c>c_*$;}
\end{array}
\right.
\end{equation}
see~\cite{Aronson-1980,Kamin-Rosenau-2004} for the one-dimensional case, and~\cite{Audrito-Vazquez-Preprint} for higher dimensions. An analogous result for the semilinear case $m=1$ was obtained by Aronson and Weinberger in their classical paper~\cite{Aronson-Weinberger-1978} (see~\cite{Aronson-Weinberger-1975} for the case of one space dimension). 

  In the $m=1$ case, more accurate estimates on the exact location of the level sets of $u$ (the transition fronts) have been obtained (see, for example,  \cite{Bramson-1978, Bramson-1983, Gartner-1982, Hamel-Nolen-Roquejoffre-Ryzhik-2013, Lau-1985, Uchiyama-1978, Uchiyama-1985}), where a ``logarithmic shift'' phenomenon  has been proved to exist. For the $m>1$ case, however, the exact location of such level sets has not been determined (for $N\geq 2$) and this question was raised again recently in Audrito and V\'{a}zquez \cite{Audrito-Vazquez-Preprint}  as an open problem (see section 12 there).

In this paper, we give an answer to this question.
In the radial case we will give a precise description of the positivity set of the solution to \eqref{eq:main} for large times, which can be used to significantly sharpen the above result in \eqref{spreading-rough}  for the general case, revealing the precise logarithmic shift.

We now describe our main results more accurately. 
Given a radially symmetric and compactly supported nonnegative initial function $u_0\not\equiv 0$, the solution of \eqref{eq:main} is also radial and compactly supported for any positive time, and satisfies the equation
$$
u_t=(u^m)_{rr}+\frac{N-1}{r}(u^m)_r+u(1-u),
$$
where $r=|x|$  and, abusing notation, $u=u(r,t)$. Let
\begin{equation}
\label{eq:def.h}
h(t)=\inf\{r>0: u(x,t)=0\text{ if }|x|>r\}.
\end{equation}
It is easy to show that after some finite time $T$ the spatial support of $u(\cdot, t)$ for any later time is a ball of radius $h(t)$: $u(x,t)>0$ if $|x|<h(t)$, and $u(x,t)=0$ if $|x|\ge h(t)$ for all $t\ge T$.

To determine the precise long-time behavior of $h(t)$, the main ideas in our approach can be briefly explained as follows.
 From \eqref{spreading-rough} we have $\lim_{t\to\infty}\frac{h(t)}t=c_*$. Thus, close to the free boundary $r=h(t)$,
\[
\frac{N-1}{r}\approx \alpha(t):=\frac{N-1}{c_*t},
\]
and $u$ solves approximately the nonlinear reaction-diffusion-advection equation
$$
u_t=(u^m)_{rr}+\alpha(t)(u^m)_r+u(1-u).
$$
When $\alpha(t)$ is replaced by a constant $\alpha$, this problem is known to have one-dimensional wavefronts $\Phi_{c,\alpha}$ for all speeds $c\ge c(\alpha)$, where $c(\alpha)$ is a critical speed depending on $\alpha$ (and $m$, of course); see for instance the monograph~\cite{Gilding-Kersner-2004} or the paper~\cite{Malaguti-Ruggerini-2010} and the references therein. For $c=c(\alpha)$ there is a unique  wavefront $\Phi_{c(\alpha), \alpha}$ with this speed whose support is $\overline{\mathbb{R}}_-$.
One would guess that
$$
h'(t)\approx c(\alpha(t))\approx c(0)+c'(0)\frac{N-1}{c_*t}
$$
for large times.
Since $c(0)=c_*$, we expect that
$$
h(t)\approx c_*t -(N-1)c^*\log t\ \text{with }c^*=-c'(0)/c_*>0
$$
as $t\to\infty$.
Moreover,  $\Phi_{c(\alpha), \alpha}$ approaches $\Phi_{c_*}$ as $\alpha\to0^+$. Hence, we also expect the solution $u(r,t)$ to approach the profile $\Phi_{c_*}$ in a suitable moving coordinate system as $t\to\infty$. We will show that this is indeed the case.
\begin{Theorem}
\label{thm:main}
Let $u$ be the solution to~\eqref{eq:main} with a radially symmetric and compactly supported nonnegative initial function $u_0\in L^\infty(\mathbb{R}^N)\setminus\{0\}$, and let $h$ be the function describing the interface, defined by~\eqref{eq:def.h}.
Then there is a constant $r_0$ such that
\begin{equation}
\label{eq:lim-h(t)}
\lim_{t\to\infty} \big[h(t)-c_*t+(N-1) c^*\log t\big]=r_0,
\end{equation}
\begin{equation}
\label{eq:convergence.in.form}
\lim_{t\to\infty}\left\{\sup_{r\geq 0}\big|u(r,t)-\Phi_{c_*}(r-c_*t+(N-1)c^*\log t-r_0)\big|\right\}=0.
\end{equation}
\end{Theorem}

\medskip

An explicit expression for $c^*$ is given in Section~\ref{sect:convection}; see \eqref{c^*}.

\medskip

We note that when the spatial dimension  $N=1$, the coefficient of the $\log t$ term in~\eqref{eq:lim-h(t)} and~\eqref{eq:convergence.in.form} vanishes, and hence no logarithmic correction occurs. This was already suggested
by the paper~\cite{Biro-2002}, where it is shown that\footnote{See Section 6.3 below regarding a gap in the proof of this result in \cite{Biro-2002}, and how the gap can be fixed.}
\[\lim_{t\to\infty} \left\{\sup_{x\in \mathbb R}\big|u(x,t)-\Phi_{c_*}(x-c_*t-x_0)\big|\right\}=0 \mbox{
 for some $x_0\in \mathbb R$,}
\]
if the initial function in \eqref{eq:main} is  nonnegative, piecewise continuous, bounded, $u_0(x)=0$ for $x\ge \bar x$ for some $\bar x\in \mathbb{R}$, and $\liminf_{x\to-\infty}u_0(x)>0$.
Such a phenomenon is in  sharp contrast to
 the semilinear case $m=1$. Indeed, it follows from Bramson~\cite{Bramson-1978,Bramson-1983} that in space dimension one with $m=1$ and $u_0$ as  above, $$
\lim_{t\to\infty}\left\{\sup_{x\in\mathbb R}|u(x,t)-\Phi_{c_*}(x-c_*t+\frac3{c_*}\log t-r_0)|\right\}=0,
$$
where $\Phi_{c_*}$ is the (positive) wavefront for the equation with minimal speed $c_*=2$ such that $\Phi_{c_*}(0)=1/2$;  see also~\cite{Hamel-Nolen-Roquejoffre-Ryzhik-2013,Lau-1985,Uchiyama-1978}. Thus for this kind of initial functions, the level sets of $u(\cdot, t)$ move asymptotically with speed $c_*-\frac{3}{c_*t}$ instead of exactly $c_*$. The associated logarithmic shifting term $\frac{3}{c_*}\log t$ is known in the literature as the logarithmic Bramson correction term.
For higher dimensions with radially symmetric initial function $u_0$ as in Theorem \ref{thm:main}, the results in~\cite{Gartner-1982,Uchiyama-1985} indicate that Bramson's correction term for the case $m=1$ becomes $\frac{N+2}{c_*}\log t$. Theorem~\ref{thm:main} shows that such a correction is also present for the porous medium case $m>1$, though with a coefficient of the form $c^*(N-1)$. This resembles the case of the Stefan problem with a Fisher-KPP term, where the logarithmic correction term is also given by $(N-1)c^*\log t$ (with a different value for $c^*$); see~\cite{Du-Matsuzawa-Zhou-2015}, and also~\cite{Du-Matsuzawa-Zhou-2014} for $N=1$. 

Theorem \ref{thm:main} can be easily applied to understand the spreading speed of \eqref{eq:main} with initial functions which are not radially symmetric. Indeed, suppose that $u$ is the solution of \eqref{eq:main}
with initial function $u_0$ bounded, nonnegative and having nonempty compact support. Let $\underline u_0$ and $\overline u_0$ be nonnegative  radially symmetric functions in $L^\infty(\mathbb R^N)\setminus\{0\}$ with compact support such that, for some $t_0>0$,
\[
 \underline u_0(x)\leq u(x, t_0)\leq \overline u_0(x) \mbox{ in } \mathbb R^N,
\]
and let $\underline u$ and $\overline u$ be the solutions of \eqref{eq:main} with initial functions $\underline u_0$ and $\overline u_0$, respectively. Then by the comparison principle we have
$\underline u(x,t)\leq u(x,t+t_0)\leq \overline u(x, t)$ in $\mathbb R^N\times \mathbb R_+$.  
Applying Theorem~\ref{thm:main} to $\underline u$ and $\overline u$, we immediately obtain  the following sharper version of~\eqref{spreading-rough}.
\begin{Corollary}\label{cor:1} For $u_0$ as described above,
there exist $r_1, r_2\in\mathbb R$ such that  the boundary of the spatial support of the solution $u(\cdot, t)$ for all large  time $t$ is contained in the spherical shell
\[
\big\{x\in\mathbb R^N: r_1\leq |x|-c_* t+(N-1)c^* \log t\leq r_2\big\},
\]
and
\[\lim_{t\to\infty}u(x,t)=1 \mbox{ uniformly in } \big\{|x|\leq c_*t-(N-1)c\log t\big\} \text{ for any }c>c^*.
\]
\end{Corollary}

\noindent\textsc{Organization of the paper. } In Section 2 we will consider the wavefront  $\Phi(x;\alpha)$ with minimal speed $c(\alpha)$ for the corresponding one dimensional reaction-diffusion-advection equation of \eqref{eq:main}, with advection term $\alpha (u^m)_x$.
We will study the properties of $c(\alpha)$ and $\Phi(x;\alpha)$ as  functions of $\alpha\in\mathbb R$, showing in particular that they vary smoothly with $\alpha$; this lays down the foundation of our approach.
In Section 3 we will use suitable modifications of $\Phi$ to construct  subsolutions  and supersolutions  to show  that the position of the  free boundary has a logarithmic correction, with coefficient $c^*(N-1)$. To prove that as $t\to\infty$ the front converges
in the moving frame $r=c_*t-c^*(N-1)\log t$, much extra work is needed. Firstly we prove a uniform bound for the flux $(u^m)_r$ in Section 4, then
we show  in Section 5 that certain eternal solutions of \eqref{eq:main} must coincide with wavefronts. 
Finally, with the help of these results, we show in Section 6.1 that there exists a time sequence $t_n\to\infty$ such that along this sequence, 
$u$ converges to a finite shift of $\Phi_{c_*}=\Phi(\cdot; 0)$ in the above mentioned moving frame; this enables us to refine the subsolutions and supersolutions constructed in Section~3 to show in Section 6.2 that the convergence of $u$ holds for $t\to\infty$. In Section~6.3, we discuss a gap in \cite{Biro-2002} and also give a version of Theorem 1.1 in dimension one, without requiring $u_0$ to be symmetric.

\section{Wavefronts for the convection problem}
\label{sect:convection}
\setcounter{equation}{0}

The aim of this section is to study the dependence on $\alpha$ of traveling wave solutions to the one dimensional convection problem
$$
u_t=(u^m)_{xx}+\alpha (u^m)_x+u(1-u).
$$
From Theorem 3.1 in~\cite{Malaguti-Ruggerini-2010} we know that for any constant $\alpha\in\mathbb R$ there exists a minimum speed $c(\alpha)$ corresponding to which there is a unique wavefront $\Phi(x;\alpha)$ with support $(-\infty, 0]$.
It will be convenient to work with the pressure variable
$$
\phi(x;\alpha):=\frac{m}{m-1}\Phi^{m-1}(x;\alpha),
$$
which satisfies $\phi\in C(\mathbb{R})\cap C^\infty(\mathbb{R}_-)$, and
\begin{equation}
\label{eq:TW.convection}
\left\{
\begin{array}{l}
-c(\alpha)\phi'-(m-1)\phi\phi''-(\phi')^2-\alpha(m-1)\phi\phi'=f(\phi),\quad  \phi>0 \quad\text{in }\mathbb{R}_-,\\[6pt]
\phi=0\quad\text{in }\mathbb{R}_+,\qquad\phi(-\infty)=\frac{m}{m-1},\quad \big(\phi^{\frac{m}{m-1}}\big)'(0^-)=0,
\end{array}
\right.
\end{equation}
where
$$
f(\phi):=(m-1)\phi\Big[1-\big(\frac{m-1}{m}\phi\big)^{\frac{1}{m-1}}\Big].
$$

To simplify notations,   here and in what follows we will often omit the $\alpha$ dependence for functions of the form $\eta(s;\alpha)$, except when confusion may arise. We will moreover use $\eta'(s;\alpha)$ to denote its derivative with respect to $s$,
and  $\partial_\alpha\eta (s;\alpha)$ or $\eta_\alpha(s;\alpha)$ to denote its derivative with respect to $\alpha$.

The solutions to the auxiliary problem~\eqref{eq:TW.convection} will  be used to construct suitable sub- and supersolutions to problem~\eqref{eq:main}. Many results on problem~\eqref{eq:TW.convection} can be found  in~\cite{Malaguti-Marcelli-Matucci-2010, Malaguti-Ruggerini-2010}. Here we will obtain further  results, both for the velocity $c(\alpha)$ and for the profile $\phi$ as functions of $\alpha$.

\begin{Lemma}\label{lem:phi'}
The solution $\phi$ to~\eqref{eq:TW.convection} satisfies $\phi'(0;\alpha):=\lim_{x\to0^-}\phi'(x;\alpha)=-c(\alpha)$.
\end{Lemma}
\begin{proof} This follows from~\cite{Malaguti-Ruggerini-2010} with $h(u)=-\alpha mu^{m-1}$ and $f(u)=u-u^2$.
More precisely, if we denote $c(\alpha)$ by $c^*$ and $\phi(\xi;\alpha)$ by $\varphi(\xi-\beta^-)$ to match the notations in that paper,
then, from the proofs of Theorems 3.1 and 3.2 there, we have
\[
\varphi'(\xi)=\frac{z_{c^*}(\varphi)}{m\varphi^{m-1}}
\]
with $z_{c^*}(\eta)$ satisfying
\[
z_{c^*}(\eta)<0 \mbox{ for } \eta\in (0,1),\; z_{c^*}(0)=z_{c^*}(1)=0,\; z_{c^*}'(0)=-c^*.
\]
We thus obtain
\[
\lim_{x\to 0^-}\phi'(x;\alpha)=\lim_{\xi\to \beta^-}m\varphi^{m-2}(\xi)\varphi'(\xi)=\lim_{\varphi\to 0^+}\frac {z_{c^*}(\varphi)}{\varphi}=z_{c^*}'(0)=-c^*=-c(\alpha).
\]

\end{proof}

\begin{Lemma}
\label{lem:second.derivative.phi}
Let $\phi$ be the solution to~\eqref{eq:TW.convection}. Then
$\phi''(0;\alpha):=\lim_{x\to0^-}\phi''(x;\alpha)=\frac{m-1}m(\alpha c(\alpha)-1)$.
\end{Lemma}
\begin{proof}
We first assume that $\lim_{x\to0^-}\phi''(x)=\beta\in\mathbb{R}\cup\{\pm\infty\}$ exists. We may rewrite the first equation
in \eqref{eq:TW.convection} in the form
$$
\phi''=\frac{-c(\alpha)\phi'-(\phi')^2}{(m-1)\phi}-\alpha\phi'-1+\Big(\frac{m-1}{m}\Big)^{\frac{1}{m-1}}\phi^{\frac1{m-1}}.
$$
Let $x\to0^-$. Then $\phi\to0$, and by L'H\^opital's rule and Lemma \ref{lem:phi'} we get
$$
\beta=-\frac{\beta}{m-1}+\alpha c(\alpha)-1, \mbox{ and hence } \beta=\frac{m-1}m(\alpha c(\alpha)-1).
$$

If the limit does not exist, then $\limsup_{x\to 0^-}\phi''(x)>\liminf_{x\to 0^-}\phi''(x)$. Hence there exists a sequence of local maxima of $\phi''(x)$, denoted by $\{y_n\}$, such that
\[
y_n\to 0^-,\; \phi''(y_n)\to\limsup_{x\to 0^-}\phi''(x),\; \phi'''(y_n)=0.
\]
Similarly there exists a sequence $\{z_n\}$ satisfying
\[
z_n\to 0^-,\; \phi''(z_n)\to\liminf_{x\to 0^-}\phi''(x),\; \phi'''(z_n)=0.
\]

Differentiating the equation in~\eqref{eq:TW.convection}, we have
\begin{equation}
\label{eq:differentiated.equation}
\begin{array}{l}
-c(\alpha)\phi''-2\phi'\phi''-(m-1)\phi\phi'''-(m-1)\phi'\phi''-\alpha(m-1)(\phi')^2-\alpha(m-1)\phi\phi''
\\[5pt] \qquad\qquad
=(m-1)\phi'-m\Big(\frac{m-1}m\Big)^{\frac1{m-1}}\phi^{\frac{1}{m-1}}\phi',
\end{array}
\end{equation}
and hence, for $x_n\in \{y_n, z_n\}$,
$$
\begin{array}{l}
-c(\alpha)\phi''(x_n)-(m+1)(-c(\alpha)+o(1))\phi''(x_n)-\alpha (m-1) c^2(\alpha)(1+o(1))-o(1)\phi''(x_n)\\[10pt]
\qquad\qquad
=-(m-1)c(\alpha)(1+o(1))+o(1),
\end{array}
$$
that is,
$$
\phi''(x_n)=\frac{m-1}m(\alpha c(\alpha)-1)(1+o(1)).
$$
It follows that
\[
\lim_{n\to\infty} \phi''(y_n)=\lim_{n\to\infty} \phi''(z_n)=\frac{m-1}m(\alpha c(\alpha)-1).
\]
Hence  $\limsup_{x\to 0^-}\phi''(x)=\liminf_{x\to 0^-}\phi''(x)$.
This contradiction completes the proof.
\end{proof}

\begin{Lemma} The function $c(\alpha)$ is a Lipschitz continuous function on $\alpha$, and
\begin{equation*}
\label{eq:lipschitz.monotone}
0\le c(\alpha_1)-c(\alpha_2)\le m(\alpha_2-\alpha_1)\quad\text{for any } \; \alpha_2>\alpha_1.
\end{equation*}
\end{Lemma}

\begin{proof}This follows easily from the proof of Proposition 3.1 in~\cite{Malaguti-Ruggerini-2010} and the fact that
$0\leq \alpha_2m u^{m-1}-\alpha_1 m u^{m-1}\leq m(\alpha_2-\alpha_1)$ for $u\in [0,1]$ and $\alpha_1<\alpha_2$.
\end{proof}

To obtain further properties of the function $c(\alpha)$, it is convenient to work on the derivative of $\phi(\cdot;\alpha)$ at the level $q$, which we denote by
$$
p(q;\alpha):=\phi'(\phi^{-1}(q;\alpha);\alpha),\quad q\in\big(0,\frac{m}{m-1}\big),\ \alpha\ge0.
$$
Let us remark that the profiles $\phi$ are strictly monotone and smooth in $\mathbb{R}_-$; see~\cite{Gilding-Kersner-2004}. Hence $p(q;\alpha)$ is well defined.

We have, using~\eqref{eq:TW.convection},
\begin{equation}
\label{eq:derivative.p}
p'(q;\alpha)=\frac{\phi''(\phi^{-1}(q;\alpha);\alpha)}{\phi'(\phi^{-1}(q;\alpha);\alpha)}
=-\frac{c(\alpha)}{(m-1)q}-\frac{p(q;\alpha)}{(m-1)q}-\alpha-\frac{f(q)}{(m-1)q\, p(q;\alpha)}.
\end{equation}
Clearly
\begin{equation}\label{p-0}
\begin{aligned}
&p(0;\alpha)=\phi'(0;\alpha)=-c(\alpha),\; p'(0;\alpha)=\phi''(0;\alpha)=\frac{m-1} m (-\alpha +\frac{1}{c(\alpha)}),
\\
&p(\frac{m}{m-1};\alpha):=\lim_{q\to (\frac{m}{m-1})^-} p(q;\alpha)=\lim_{x\to-\infty}\phi'(x;\alpha)=0.
\end{aligned}
\end{equation}
\begin{Lemma}
\label{lem:derivative.p.alpha}
For any $\alpha\in\mathbb R$,
$$
p'\Big(\frac{m}{m-1};\alpha\Big):=\lim_{q\to (\frac{m}{m-1})^-}p'(q;\alpha)=\gamma(\alpha) \mbox{ with }
$$
$$
\gamma(\alpha):=\frac{c(\alpha)+m\alpha+\sqrt{(c(\alpha)
+m\alpha)^2+4m}}{2m}.
$$
\end{Lemma}
\begin{proof}
If  $\lim_{q\to\frac{m}{m-1}}p'(q;\alpha)=\beta$ exists, then from~\eqref{eq:derivative.p} we get, using l'H\^opital's rule and also the fact $f'\big(\frac{m}{m-1}\big)=-1$, that $\beta$ is a solution to
\begin{equation*}
\label{eq:satisfied.by.gamma}
\beta=-\frac{c(\alpha)}m-\alpha+\frac{1}{m\beta}.
\end{equation*}
This equation has two roots, one positive, $\beta_+$, and one negative, $\beta_-$. Since $p(q;\alpha)<0$ for $r\in (0, \frac{m}{m-1})$ and $p(\frac m{m-1};\alpha)=0$, the limit must be given by the positive one, $\beta_+=\gamma(\alpha)$.

Let us now prove that the limit exists. It is convenient to work with the function
\[
\psi(x):=p'(\phi(x);\alpha)=\frac{\phi''(x)}{\phi'(x)},\; x\in\mathbb{R}_-.
\]
It suffices to show that $\lim_{x\to-\infty}\psi(x)$ exists. Otherwise
\[
\limsup_{x\to-\infty}\psi(x)>\liminf_{x\to-\infty}\psi(x).
\]
So, as in the proof of Lemma \ref{lem:second.derivative.phi}, we can find sequences $\{y_n\}$ and $\{z_n\}$ such that
\[
\lim_{n\to\infty} y_n=\lim_{n\to\infty} z_n=-\infty,\; \psi'(y_n)=\psi'(z_n)=0,\]
\[
 \psi(y_n)\to \limsup_{x\to-\infty}\psi(x),\; \psi(z_n)\to \liminf_{x\to\infty}\psi(x).
\]

Equation~\eqref{eq:differentiated.equation} can be rewritten as
\begin{equation}
\label{eq:phi-2}
(m-1)\phi\frac{\phi'''}{\phi'}+\big[c(\alpha)+(m+1)\phi'+\alpha(m-1)\phi\big]\frac{\phi''}{\phi'}\smallskip
=m-1-m\left(\frac{m-1}m\phi\right)^{\frac1{m-1}}+\alpha\phi'.
\end{equation}
Since
\[
\frac{\phi'''}{\phi'}=\psi'+\psi^2, \mbox{ and near } x=-\infty, \; \phi=\frac{m}{m-1}+o(1),\; \phi'=o(1),
\]
identity~\eqref{eq:phi-2} gives, with $x_n\in\{y_n, z_n\}$,
\[
[m+o(1)]\psi^2(x_n)+[c(\alpha)+\alpha m+o(1)]\psi(x_n)=-1+o(1).
\]
Therefore if we can show that $\psi(x)$ does not change sign for all large negative $x$, then it follows from the above identity that $\lim_{n\to\infty}\psi(y_n)=\lim_{n\to\infty}\psi(z_n)$, a contradiction.

To complete the proof, it remains to show that $\psi(x)$ does not change sign for all large negative $x$.
Denoting $v(x):=\phi'(x)$, and using $\phi(x)=\frac{m}{m-1}+o(1)$, $\phi'(x)=o(1)$ for large negative $x$, we may
now rewrite \eqref{eq:phi-2} as
\[
[m+o(1)]v''+[c(\alpha)+\alpha m+o(1)]v'=[1+o(1)]v \mbox{ for large negative $x$}.
\]
Hence $v$ can not have a negative local minimum for such $x$. Since $v<0$ in $\mathbb{R}_-$, this implies that $v$ is monotone in $(-\infty, -M]$ for some large $M>0$. It follows that $\phi''(x)=v'(x)$ does not change sign for $x\leq -M$ and hence $\psi(x)=\phi''(x)/\phi'(x)$ does not change sign for $x\leq M$.
\end{proof}

\begin{Lemma}
The function $\alpha\to c(\alpha)$ belongs to $C^2(\mathbb R)$, $\partial_\alpha p(q;\alpha)$ exists and $\alpha\to\partial_\alpha p(\cdot;\alpha)$ is continuous from $\mathbb R$ to $C([0,\frac{m}{m-1}])$. Moreover, $0>c'(\alpha)>-m$ for $\alpha\in\mathbb R$.
\end{Lemma}

\begin{proof} For arbitrarily fixed $\alpha\in\mathbb R$ and $q\in (0,\frac m{m-1})$, and $k\not=0$,
we define the incremental quotients
\[
 \hat p_k(q):=\frac{p(q; \alpha+k)-p(q;\alpha)}{k},\;
\hat c_k:=\frac{c(\alpha+k)-c(\alpha)}{k}.
\]
 A simple computation yields
$$
\hat p_k'(q)+\frac1{(m-1)q} \Big[1-\frac{f(q)}{p(q;\alpha+k)p(q;\alpha)}\Big] \hat p_k(q)=-\frac{\hat c_k}{(m-1)q}-1.
$$
Let us denote
\[
a(q):=-\frac1{(m-1)q}\Big[1-\frac{f(q)}{p(q;\alpha+k)p(q; \alpha)}\Big] \mbox{ and } A(q):=\int_q^1 a(s)\,ds.
\]
Then
$$
\big(e^{A(q)}\hat p_k(q)\big)'=e^{A(q)}\big[-\frac{\hat c_k}{(m-1)q}-1\big],
$$
which, after integration over $[\xi_1,\xi_2]\subset (0,\frac{m}{m-1})$, yields
\begin{equation}
\label{eq:incremental.quotients}
e^{A(\xi_2)}\hat p_k(\xi_2)- e^{A(\xi_1)}\hat p_k(\xi_1)=-\int_{\xi_1}^{\xi_2} e^{A(q)}\Big[\frac{\hat c_k}{(m-1)q}+1\Big]\,dq.
\end{equation}

We now examine the behavior of $e^{A(q)}$ and $\hat p_k(q)$ near $q=0$ and $q=m/(m-1)$, respectively.
Near $q=0$,
\[
f(q)=[(m-1)+o(1)]q,\; p(q;\alpha)=\phi'(0)+o(1)=-c(\alpha)+o(1),
\]
and hence
\begin{align*}
a(q)&=-\frac{1}{(m-1)q}+\frac{1+o(1)}{c(\alpha+k)c(\alpha)},
\\
A(q)&=\frac{1}{m-1}\log q+\frac{1+o(1)}{c(\alpha+k)c(\alpha)},
\\
 e^{A(q)}&=q^{\frac{1}{m-1}}[\sigma+o(1)] \mbox{ with } \sigma:=e^{\frac{1}{c(\alpha+k)c(\alpha)}},\\
\hat p_k(q)&=-\hat c_k+o(1).
\end{align*}

Near $q=m/(m-1)$,
\begin{equation}\label{f-p-m/(m-1)}
f(q)=\Big(\frac m{m-1}-q\Big)[1+o(1)],\; p(q;\alpha)=-\gamma(\alpha)\Big(\frac m{m-1}-q\Big)[1+o(1)],
\end{equation}
and hence
\begin{align*}
a(q)&=\frac 1m\Big[1-\frac{1}{\gamma(\alpha+k)\gamma(\alpha)\big(\frac{m}{m-1}-q\big)}\Big][1+o(1)],
\\
A(q)&=[1+o(1)]\log\Big(\frac m{m-1} -q\Big),\\
 e^{A(q)}& =o(1),\;\;\; \hat p_k(q)=o(1).
\end{align*}

The above analysis shows that the integral on the right hand side of \eqref{eq:incremental.quotients} is convergent when $\xi_1=0$ and $\xi_2=m/(m-1)$, and the left hand side of this equation converges to 0 as $\xi_1\to 0^+$ and $\xi_2\to (\frac m{m-1})^+$.
We thus obtain
\[
0=-\int_{0}^{\frac{m}{m-1}} e^{A(q)}\Big[\frac{\hat c_k}{(m-1)q}+1\Big]\,dq
\]
and
\[
\hat c_k=-\frac{(m-1)\int_0^{\frac{m}{m-1}}e^{A(q)}dq}{ \int_0^{\frac{m}{m-1}}e^{A(q)}q^{-1}dq}.
\]

Clearly
$$
e^{A(q)}=q^{\frac1{m-1}}\textrm{exp}\left(\frac{1}{m-1}\int_q^1\frac{f(s)}{s\,p(s; \alpha+k)p(s;\alpha)}\,ds\right).
$$
Since $c(\alpha)$ and $\gamma(\alpha)$ are continuous functions of $\alpha$, by \eqref{p-0}, Lemma \ref{lem:derivative.p.alpha}
and a standard  analysis  to the solution $p(q;\alpha)$ of the singular ODE \eqref{eq:derivative.p}, we see that $\alpha\to p(\cdot;\alpha)$
is continuous from $\mathbb R$ to $C^1([0,\frac{m}{m-1}]).$\footnote{Actually, for any given $L>0$, we could extend $p(q;\alpha)$ for $q$ to a small neighorhood of the closed interval $[0,\frac{m}{m-1}]$, say $(-\epsilon, \frac{m}{m-1}+\epsilon)$ with some small $\epsilon>0$, such that $p(q;\alpha)$ and $p'(q;\alpha)$ are jointly continuous for $(q,\alpha)\in
 (-\epsilon, \frac{m}{m-1}+\epsilon)\times [-L, L]$.} It follows that
\begin{equation}
\label{eq:definition.Psi}
e^{A(q)}\to \Psi(q;\alpha):=q^{\frac 1{m-1}}\textrm{exp}\left(\frac{1}{m-1}{\int_q^1\frac{f(s)}{s\,p^2(s; \alpha)}\,ds}\right)
\end{equation}
as $k\to 0$ uniformly for $q\in [0,\frac{m}{m-1}]$.
We may now
let $k\to 0$ in the above expression for $\hat c_k$ to obtain
\begin{equation*}\label{c'}
c'(\alpha)=-\frac{(m-1)\int_0^{\frac{m}{m-1}}\Psi(q;\alpha)dq
}{\int_0^{\frac{m}{m-1}}q^{-1}
\Psi(q;\alpha)dq} .
\end{equation*}
The continuity of $c'(\alpha)$ with respect to $\alpha\in\mathbb R$ follows from the above formula and the continuity of $\alpha\to p(\cdot;\alpha)$ in the $C^1([0,\frac{m}{m-1}])$ norm. It is also easily seen from the above formula that
\[
0>c'(\alpha)>-m \mbox{ for } \alpha\in\mathbb R.
\]

From \eqref{eq:incremental.quotients} we easily obtain
\[
\hat p_k(q)=-e^{-A(q)}\int_{0}^{q} e^{A(s)}\Big[\frac{\hat c_k}{(m-1)s}+1\Big]\,ds \; \mbox{ for } q\in (0,\frac{m}{m-1}).
\]
Letting $k\to 0$ we have, for $q\in (0,\frac{m}{m-1})$,
\begin{equation}\label{p-alpha}
\begin{array}{ll}
\partial_\alpha p(q;\alpha)&=\displaystyle -\frac{1}{\Psi(q;\alpha)}\int_{0}^{q} \Psi(s;\alpha)\Big[\frac{c'(\alpha)}{(m-1)s}+1\Big]\,ds \medskip\\
&=\displaystyle \frac{1}{\Psi(q;\alpha)}\int_{q}^{\frac m{m-1}} \Psi(s;\alpha)\Big[\frac{c'(\alpha)}{(m-1)s}+1\Big]\,ds.
\end{array}
\end{equation}

From the definition of $\Psi(q;\alpha)$ we find that for $q$ near 0,
\[
\Psi(q;\alpha)=[\sigma_\alpha+o(1)]q^{\frac 1{m-1}} \mbox{ with } \sigma_\alpha:=\exp\Big({\frac 1{m-1}\int_0^1\frac{f(s)}{s p^2(s;\alpha)}ds}\Big),
\]
and for $q$ near $\frac m{m-1}$,
\[
\Psi(q;\alpha)=[\tilde\sigma_\alpha+o(1)]\log\big(\frac m{m-1}-q\big),
\]
with $\tilde\sigma_\alpha:=\frac{-1}{m\gamma^2(\alpha)}(\frac m{m-1})^{\frac 1{m-1}}$.
It then follows from \eqref{p-alpha}  that
\[
\partial_\alpha p(0;\alpha):=\lim_{q\to 0^+}\partial_\alpha p(q;\alpha)=c'(\alpha),
\]
and for $q$ near $\frac m{m-1}$,
\begin{equation}\label{p_alpha-asymp}
\partial_\alpha p(q;\alpha)=\big[1+\frac{c'(\alpha)}m +o(1)\big]\big(\frac m{m-1} -q\big).
\end{equation}
Therefore $\alpha\to\partial_\alpha p(\cdot;\alpha)$ is continuous from $\mathbb R$ to $C([0,\frac{m}{m-1}])$.
Using this fact and the expression of $\Psi(q;\alpha)$ it is easy to check that $\partial_\alpha \Psi(q;\alpha)$ exists
and $\alpha\to \partial_\alpha \Psi(\cdot;\alpha)$ is continuous from $\mathbb R$ to $C([0,\frac{m}{m-1}])$.
This in turn implies, by the formula for $c'(\alpha)$ above, that $c''(\alpha)$ exists and is continuous in $\alpha$.
\end{proof}

\begin{Remark}\label{rm:1}
Define
\[
\Upsilon(x;\alpha):=\Psi(\phi(x); \alpha) \mbox{ with } \phi(x)=\phi(x;\alpha),
\]
with $\Psi$ as in~\eqref{eq:definition.Psi}.
Then by a simple change of variable calculation, we obtain
\[
\Upsilon(x;\alpha)=\phi(x)^{\frac{1}{m-1}}\exp\left(\frac{1}{m-1}\int_x^{\phi^{-1}(1)}\frac{f(\phi(y))}{\phi(y)\phi'(y)}dy\right)
\]
and
\[c'(\alpha)=-\frac{(m-1)\int_{-\infty}^0\phi'(x)\Upsilon(x;\alpha)dx}{\int_{-\infty}^0\phi(x)^{-1}\phi'(x)\Upsilon(x;\alpha)dx}.
\]
In particular, for $\alpha=0$, using
\[
\phi(x;0)=\frac{m}{m-1}\Phi_{c_*}^{m-1}(x) \mbox{ and } f(\phi)=(m-1)\phi\Big[1-\big(\frac{m-1}m \phi\big)^{\frac1{m-1}}\Big],
\]
we obtain, after a simple calculation, that the constant $c^*$ in Theorem~\ref{thm:main} is given by
\begin{equation}\label{c^*}
c^*=-\frac{c'(0)}{c_*}=\frac{1}{c_*}\,\frac{\int_{-\infty}^0(\Phi^m_{c_*})'(x)\,{\exp}\left(\frac{m-1}{m}\int_x^{x_*}\frac{1-\Phi_{c_*}(y)}{(\Phi_{c_*}^{m-1})'(y)}dy\right)dx}{\int_{-\infty}^0(\Phi_{c_*})'(x)\,{\exp}\left(\frac{m-1}{m}\int_x^{x_*}\frac{1-\Phi_{c_*}(y)}{(\Phi_{c_*}^{m-1})'(y)}dy\right)dx},
\end{equation}
with
$$
x_*\in \mathbb{R}_- \mbox{ uniquely determined by } \Phi_{c_*}(x_*)=\Big(\frac{m-1}{m}\Big)^{\frac{1}{m-1}}.
$$

\end{Remark}

The following technical result will play a crucial role in our analysis in the next section.

\begin{Lemma}
\label{lem:phi.alpha.order.1} For any $L>0$, there exists $C=C_L>0$ such that
\[
|\partial_\alpha \phi(x;\alpha)|\leq C \mbox{ for } x\leq 0,\; |\alpha|\leq L.
\]
\end{Lemma}

\begin{proof} Since
\[
p(\phi(x;\alpha); \alpha)=\phi'(x;\alpha) \mbox{ and } \phi(0;\alpha)=0,
\]
and since $q\to p'(q;\alpha)$ is continuous over a small neighborhood of $[0, \frac m{m-1}]$ (see the footnote in the previous page),
$\phi(x;\alpha)$ must be identical to the unique solution of the initial value problem
\[
\phi'=p(\phi;\alpha),\; \phi(0)=0
\]
in the range $x\leq 0$.
Since $\partial_\alpha p(q;\alpha)$ is continuous in both variables, by standard ODE theory (see, for example, Chapter I of \cite{Coppel}),
$\partial_\alpha \phi(x;\alpha)$ exists and $\psi(x;\alpha):=\partial_\alpha\phi(x;\alpha)$ is the unique solution of
\begin{equation}\label{phi_alpha}
\psi'=\xi(x;\alpha)\psi+\eta(x;\alpha),\; \psi(0)=0,
\end{equation}
 with
\[
\xi(x;\alpha):=p'(\phi(x;\alpha);\alpha),\;\eta(x;\alpha):=\partial_\alpha p(\phi(x;\alpha);\alpha).
\]

We claim that for any $L>0$, there exists $x_L<0$ such that
\begin{equation}\label{psi>0}
\psi(x;\alpha)>0 \mbox{ for } x\leq x_L,\; |\alpha|\leq L.
\end{equation}
Otherwise there exist a sequence $\alpha_n\in [-L, L]$ and a sequence $x_n\to-\infty$ such that $\psi(x_n;\alpha_n)\leq 0$
for all $n\geq 1$. By passing to a subsequence we may assume that $\alpha_n\to \alpha_0$ as $n\to \infty$.

By Lemma \ref{lem:derivative.p.alpha} and \eqref{p_alpha-asymp}, there exist $x_*<0$ large negative
and $\epsilon>0$ small such that
\begin{equation}\label{xi-eta>0}
\xi(x;\alpha)>0,\; \eta(x;\alpha)>0 \mbox{ for } x\leq x_*,\; |\alpha-\alpha_0|\leq 2\epsilon.
\end{equation}
Choose $n_0$ large such that $x_{n_0}\leq x_*$ and $|\alpha_n-\alpha_0|<\epsilon$ for all $n\geq n_0$. Then by \eqref{phi_alpha} we easily see that $\psi(x;\alpha_{n_0})<0$
for $x<x_{n_0}$. We now fix $x_0<x_{n_0}$ and obtain, by shrinking $\epsilon>0$ if necessary, that
$\psi(x_0;\alpha)<0$ for $|\alpha-\alpha_{n_0}|\leq \epsilon$. Using \eqref{phi_alpha} again we deduce
\[
\psi(x;\alpha)<0 \mbox{ for } x\leq x_0,\; \alpha\in [\alpha_{n_0}-\epsilon, \alpha_{n_0}+\epsilon],
\]
which implies, in view of $\psi(x;\alpha)=\partial_\alpha \phi(x;\alpha)$,
\begin{equation}\label{alpha_0-diff}
\phi(x;\alpha_{n_0}-\epsilon)>\phi(x; \alpha_{n_0}+\epsilon) \mbox{ for } x\leq x_0.
\end{equation}

On the other hand, from \eqref{f-p-m/(m-1)} we obtain
\[
\phi'(x;\alpha)=-[\gamma(\alpha)+o(1)]\Big(\frac m{m-1}-\phi(x;\alpha)\Big) \mbox{ for $x$ near $-\infty$}.
\]
It follows that for any $\delta>0$ small, there exists $x^*=x^*_{\delta,\alpha}<0$ large negative such that
\[
0< \frac m{m-1} -\phi(x;\alpha)\leq \Big[\frac m{m-1} -\phi(x^*;\alpha)\Big]e^{(\gamma(\alpha)-\delta)(x-x^*)} \mbox{ for } x<x^*,
\]
and
\[
\frac m{m-1} -\phi(x;\alpha)\geq \Big[\frac m{m-1} -\phi(x^*;\alpha)\Big]e^{(\gamma(\alpha)+\delta)(x-x^*)} \mbox{ for } x<x^*,
\]

Since $0>c'(\alpha)>-m$, the function $c(\alpha)+m\alpha$ is strictly increasing and so, by the formula for $\gamma$, we see that it is strictly increasing. Therefore, for small enough $\delta>0$, we have
\[
\gamma(\alpha_{n_0}+\epsilon)-\delta>\gamma(\alpha_{n_0}-\epsilon)+\delta.
\]
We may then apply the above inequalities for $\frac m{m-1} -\phi(x;\alpha)$ to obtain
\[
\limsup_{x\to-\infty}\frac{\frac m{m-1} -\phi(x;\alpha_{n_0}+\epsilon)}{\frac m{m-1} -\phi(x;\alpha_{n_0}-\epsilon)}\leq 0.
\]
Therefore
\[
\phi(x; \alpha_{n_0}-\epsilon)<\phi(s; \alpha_{n_0}+\epsilon) \mbox{ for  all large negative $x$}.
\]
This is a contradiction to \eqref{alpha_0-diff} and thus \eqref{psi>0} is  proved.

Using \eqref{xi-eta>0} and a finite covering argument, we easily see that for any given $L>0$ there exists $m_L<0$
large negative such that
\[
\xi(x;\alpha)>0,\; \eta(x;\alpha)>0 \mbox{ for } x\leq m_L,\; |\alpha|\leq L.
\]
We may assume that $m_L\leq x_L$ so that we also have
$\psi(x; \alpha)>0$ in this range. We thus have
\[
\Big(\psi(x;\alpha) e^{ \int_x^{m_L}\xi(s;\alpha)ds}\Big)'=e^{\int_x^{m_L}\xi(s;\alpha)ds}\eta(x;\alpha)>0 \mbox{ for } x<m_L,\; |\alpha|\leq L.
\]
Since $x\to e^{\int_x^{m_L}\xi(s;\alpha)ds}$ is decreasing for $x\in (-\infty, m_L]$, the above inequality implies that
$x\to \psi(x;\alpha)$ is increasing for $x\in (-\infty, m_L]$ and every fixed $\alpha\in [-L, L]$. In view of $\psi(x;\alpha)>0$ in this range,
we can conclude that
\[
\sup_{x\leq 0,\; |\alpha|\leq L} |\psi(x;\alpha)|=C_L:=\max_{x\in [m_L, 0], \alpha\in [-L,L]}|\psi(x;\alpha)|.
\]
This completes the proof.
\end{proof}


\section{Logarithmic shift}
\label{sect:shift}
\setcounter{equation}{0}

Let $h(t)$ be defined by \eqref{eq:def.h}, and $c^*$ be given by \eqref{c^*}.
In this section we prove the following logarithmic shift result for $h(t)$.

\begin{Theorem}\label{thm:shift}
There exist positive constants $T$ and $C$ such that, for $t\geq T$,
\[
c_*t-(N-1)c^*\log t-C\leq h(t)\le c_* t-(N-1)c^*\log t+C.
\]
\end{Theorem}

This theorem will be proved by a sequence of lemmas.
We start by obtaining an estimate on how the solution $u(r,t)$ approaches 1 in sets of the form $0\leq r\le ct$, with $c>0$ small.

\begin{Lemma} $\lim_{t\to\infty}h(t)/t=c_*$, and
there exist $\hat c\in(0,c_*)$, $\delta\in(0,1)$ and $M,T_*>0$
such that
\begin{align}
u(r,t)\leq 1+M e^{-\delta t} & \mbox{ for all $r\geq 0$ and $t\geq T_*,$}&\label{u<}\\
u(r,t)\ge 1-Me^{-\delta t} & \mbox{ for all $r\in[0,\hat c\,t]$ and $t\ge T_*$.}&\label{u>}
\end{align}

\end{Lemma}

\begin{proof} The fact that  $\lim_{t\to\infty}h(t)/t=c_*$ is a direct consequence of \eqref{spreading-rough} and the proof of Theorem 2.6 in~\cite{Audrito-Vazquez-Preprint}.

Let $\overline u$ be the unique solution of the ODE problem
\[
\overline u'=\overline u-\overline u^2,\; \overline u(0)=\|u_0\|_\infty.
\]
Then it is well known (and easily seen) that
\[
\overline u(t)\leq 1+M e^{-\delta t} \mbox{ for all $t>0$ and some $M,\, \delta>0$.}
\]
By the comparison principle we have $u(r,t)\leq \overline u(t)$ for all $r\geq 0$ and $t>0$, which clearly implies \eqref{u<}.

We next prove \eqref{u>}. Setting $v=u^m$ we have
\[
\frac{1}{m}v^{\frac 1m -1}  v_t-\Delta v=v^{\frac 1m}(1-v^{\frac 1m}).
\]
Since
\[
\lim_{v\to 1} \frac{v^{\frac 1m}(1-v^{\frac 1m})}{1-v}=\frac 1m,
\]
there exists $\epsilon>0$ small so that
\[
v^{\frac 1m}(1-v^{\frac 1m})>\frac{1}{2m}(1-v) \mbox{ for } v\in[1- \epsilon, 1).
\]

Fix $\tilde c\in (0, c_*)$. By \eqref{spreading-rough} we have
\[
v(r,t)\to 1 \mbox{ uniformly in $\{0\leq r\leq \tilde c t\}$ as $t\to\infty$}.
\]
It follows that for all large $t$, say $t\geq T_1$,
\begin{equation}\label{v-estimates}
v\geq 1-\epsilon,\; \frac{1}{m}v^{\frac 1m -1}\leq \frac{2}{m}\;
\mbox{ for } 0\leq r\leq \tilde c t.
\end{equation}

For any fixed $T>0$ we now consider an auxiliary problem
\[
\left\{
\begin{aligned}
&\psi_t-\Delta \psi=\frac{1}{2m}(1-\psi) & \mbox{ for }& |x|<\tilde c T,\; t\in (0, T],\\
&\psi=1-\epsilon & \mbox{ for } &|x|=\tilde c T,\; t\in (0, T],\\
&\psi=1-\epsilon & \mbox{ for }& |x|\leq \tilde c T,\; t=0.
\end{aligned}
\right.
\]
One easily sees that $\psi\equiv 1-\epsilon$ is a  subsolution and $\psi\equiv 1$ is a supersolution of the above problem.
It follows that the unique solution $\psi(x,t)$ of this problem satisfies
\begin{equation}\label{psi-estimates}
\psi_t\geq 0 \mbox{ and } 1-\epsilon\leq \psi<1 \mbox{ for } |x|\leq \tilde c T,\; t\in [0, T].
\end{equation}
Moreover, by the proof of Lemma 2.6 in \cite{Du-Polacik}, there exist constants $c_1\in (0,\tilde c)$, $c_2\in (0,1)$, $\delta_0>0$, $T_2\geq T_1$ and $\hat M>0$ such that
\[
\psi(x, c_2 T)\geq 1-\hat M e^{-\delta_0 T}  \mbox{ for } |x|\leq c_1 T,\; T\geq T_2,
\]
or, equivalently,
\begin{equation}
\label{psi-estimate}
\psi(x, T)\geq 1 -\hat M e^{-\tilde \delta_0 T}  \mbox{ for } |x|\leq \tilde c_1 T,\; T\geq \tilde T_2,
\end{equation}
with $\tilde \delta_0:=\delta_0/c_2,\; \tilde c_1:=c_1/c_2,\; \tilde T_2:=T_2/c_2>T_2\geq T_1$.

For arbitrarily fixed $T\geq \tilde T_2$  we define
$
V(x, t):=v(|x|, t+T)$ and
\[
W(x,t):=\psi(x, \frac m 2 t) \mbox{ for } (x,t)\in \Omega_T:=\left\{(x,t)\in\mathbb{R}^{N+1}: |x|<\tilde cT,\; t\in (0, \frac 2m T]\right\}.
\]
Then by \eqref{v-estimates}, \eqref{psi-estimates} and the choices of $\epsilon$ and $T_1$, we obtain
\[
(\frac 1m V^{\frac 1m -1})V_t-\Delta V=V^{\frac 1m}(1-V^{\frac 1m})  \mbox{ in } \Omega_T,\; V\geq 1-\epsilon \mbox{ on } \partial_p\Omega_T,
\]
and
\[
(\frac 1m V^{\frac 1m -1})W_t-\Delta W\leq W^{\frac 1m}(1-W^{\frac 1m}) \mbox{ in } \Omega_T,\; W= 1-\epsilon \mbox{ on } \partial_p\Omega_T,
\]
where $\partial_p\Omega_T$ denotes the parabolic boundary of $\Omega_T$. The comparison principle then yields
\[
V(x,t)\geq W(x,t) \mbox{ in } \Omega_T.
\]
In particular, making use of \eqref{psi-estimate}, we obtain
\[
V(x,\frac 2m T)\geq W(x,\frac 2m T)=\psi(x, T)\geq 1-\hat M e^{-\tilde \delta_0 T} \mbox{ for } |x|\leq \tilde c_1  T.
\]
Therefore for any $T\geq \tilde T_2$ we have
\[
v\big(r, (1+\frac 2m)T\big)\geq 1-\hat M e^{- \delta (1+\frac 2m)T} \mbox{ for } 0\leq r\leq \hat c (1+\frac2m )T
\]
with $ \delta:= \tilde \delta_0/(1+\frac 2m),\; \hat c:=\tilde c_1/(1+\frac2m)$. It follows that
\[
u^m(r,t)\geq 1-\hat M e^{- \delta t} \mbox{ for } 0\leq r\leq \hat c t,\; t\geq (1+\frac 2m)\tilde T_2,
\]
from which \eqref{u>}  easily follows.
\end{proof}

 We note that all the functions $w$ that we will use as sub- and supersolutions from now on satisfy that $(w^m)_r=0$ at their free boundary. Hence, to check that they are indeed subsolutions (respectively supersolutions), it is enough to check that
$$
\mathcal{L}w:=w_t-(w^m)_{rr}-\frac{N-1}r (w^m)_r- w(1-w)\leq 0\;\; (\mbox{respectively $\geq 0$})
$$
 in the positivity set of $w$, plus the correct ordering at the parabolic boundary.

\begin{Lemma}
There exist $M>(N-1)c^*$, $T>0$, and $\delta\in(0,1)$
such that
\[
\mbox{$h(t)\ge c_*t-M\log t$ for $t\geq T$.}
\]
 Moreover,  there exists $\bar M\ge M$ and $\bar T\geq T$
such that
\[
\mbox{$u(r,t)\ge 1-t^{-2}$ for $r=c_* t-\bar M\log t$ and $t\geq \bar T$.}
\]
\end{Lemma}

\begin{proof}
We will perform comparison in the range $r\ge \hat ct$, $t\ge T$, with $\hat c$ as in the previous lemma and some large $T$ to be determined. As subsolution we will use
$$
w(r,t)=g(t)\Phi_{c_*}(r-c_*t+M\log t),\quad g(t)=1-Me^{-\delta t},
$$
with $M$ and $\delta$ chosen suitably.

By the previous lemma, we can find $\delta_0\in (0,1)$, $M_0>(N-1)c^*$, and $T_0>0$ so that~\eqref{u>} holds for $\delta\in (0, \delta_0]$, $M\geq M_0$ and $T\geq T_0$.
Since $\Phi_{c_*}<1$, we have
\[
\mbox{$w(r,t)\le u(r,t)$ at $r=\hat ct$ for $t\geq T_0$ and $\delta\in (0, \delta_0]$, $M\geq M_0$.}
\]

We show next that, with
\[
M\geq M_1:=\max\{M_0, 2m(N-1)/\hat c\} \mbox{ and } \delta:=\min\left\{\delta_0, \min\{\frac{m-1}{4}, \frac 12\}\right\},
\]
 there exists $\epsilon>0$ independent of such  $M$ so that
\begin{equation}\label{Lw<0}
\mathcal{L}w\leq 0 \mbox{ for } r\in(\hat c t,\geq c_*t-M\log t),\; t\geq \frac 1\delta\log\big(\frac M\epsilon\big).
\end{equation}
Indeed, for $r\in (\hat ct,c_*t-M\log t)$,
$$
\begin{aligned}
\mathcal{L} w=& M\delta e^{-\delta t} \Phi_{c_*}-g(t)\Big(c_*-\frac{M}t\Big)\Phi_{c_*}'-g^m(t)(\Phi_{c_*}^m)''\\
&-\frac{N-1}{r}g^m(t)(\Phi_{c_*}^m)'-g(t)\Phi_{c_*}
(1-g(t)\Phi_{c_*}),
\end{aligned}
$$
where we have omitted the argument ``$r-c_*t+M\log t$'' in $\Phi_{c_*}$ and its derivatives to simplify the presentation.
Substituting
$$
(\Phi^m_{c_*})''=-c_*\Phi_{c_*}'-\Phi_{c_*}(1-\Phi_{c_*})
$$
into the above expression we get,
$$
\begin{aligned}
\mathcal{L} w =& \Phi_{c_*}\Big[\underbrace{M\delta e^{-\delta t} +[g^{m-1}(t)-1]g(t)(1-\Phi_{c_*})-Mg(t)e^{-\delta t}\Phi_{c_*}}_{\mathcal{A}_1}\Big]
\\
&+g(t)\Phi'_{c_*}\Big[\underbrace{c_*[g^{m-1}(t)-1]+\frac{M}{t}-\frac{N-1}{r}mg^{m-1}(t)\Phi_{c_*}^{m-1}}_{\mathcal{A}_2}\Big].
\end{aligned}
$$
Using
\[
 \lim_{x\to1}\frac{x^{m-1}-1}{x-1}=m-1
\]
we can find $\epsilon\in (0, 1/2)$ such that
\[
g(t)\geq 1-\epsilon>\frac 12,\; -2(m-1)Me^{-\delta t}\leq g^{m-1}(t)-1\leq -\frac{m-1}{2}Me^{-\delta t}
\]
 when $ Me^{-\delta t}\leq \epsilon$,  or equivalently, when $t\geq \frac 1\delta\log\big(\frac M\epsilon\big)$. For later use, we also
assume that $\epsilon>0$ has been chosen small enough such that
\[
\frac 12-c_*2(m-1) t e^{-\delta t}>0 \mbox{ for } t\geq \frac 1\delta\log\big(\frac {M_1}\epsilon\big).
\]

Therefore, for such $t$,
\[
\mathcal A_1\leq M e^{-\delta t}\left[\delta-\Big(\frac{m-1}{4}(1-\Phi_{c_*})+\frac 12 \Phi_{c_*}\Big)\right]\leq Me^{-\delta t}
(\delta-\min\big\{\frac{m-1}{4}, \frac 12\big\})\leq 0.
\]
For $t\geq \frac 1\delta\log\big(\frac {M}\epsilon\big)$ and $r\geq \hat c t$, we additionally have, in view of $M\geq M_1$ and $g, \Phi_{c_*}\in (0,1)$,
$$
\begin{aligned}
\mathcal{A}_2 \geq &\; -c_*2(m-1)Me^{-\delta t}+\frac{M}{t}-\frac{m(N-1)}{\hat ct}\\
=&\;\frac Mt \Big(1-\frac{m(N-1)}{\hat c M}-c_*2(m-1) t e^{-\delta t}\Big)\\
\geq &\;\frac Mt\Big(\frac 12-c_*2(m-1) t e^{-\delta t}\Big)
>0.
\end{aligned}
$$
We thus have, for $t\geq \frac 1\delta\log\big(\frac {M}\epsilon\big)$, $M\geq M_1$ and $r\geq \hat c t$,
\[
\mathcal{L}w=\Phi_{c_*}\mathcal{A}_1+g(t)\Phi'_{c_*}\mathcal A_2\leq 0.
\]
This proves \eqref{Lw<0}.

We now take $T\geq T_0$ large enough so that
\[
T\geq \frac 1\delta\log \big(\frac{(c_*-\hat c)T}{\epsilon \log T}\big) \;\mbox{ and }\;  \frac{(c_*-\hat c)T}{\log T}\geq M_1.
\]
Then we take
\[
M:=\frac{(c_*-\hat c)T}{\log T}.
\]
Clearly $M\geq M_1$, and $t\geq T$ implies $t\geq \frac 1\delta\log\big(\frac {M}\epsilon\big)$.
Therefore for $t\geq T\geq T_0$,
\[
\mathcal{L}w\leq 0 \mbox{ for } r\geq \hat c t,\;  \mbox{ and } w(r,t)\leq u(r,t) \mbox{ for } r=\hat c t.
\]
We also have
\[
r-c_*T-M\log T\geq 0 \mbox{ for } r\geq \hat c T,
\]
and hence $w(t, T)=0\leq u(r, T)$ for $r\geq \hat c T$.

We may now apply the comparison principle over the region $\Omega:=\{(r,t): r\geq \hat c t,\; t\geq T\}$ to conclude that
\[
u(r,t)\geq w(r,t) \mbox{ in } \Omega,
\]
 which implies in particular that $h(t)\ge c_*t -M \log t$ for $t\geq T$.

It is well known that
\[
\Phi_{c_*}(x)\geq 1-e^{\delta_* x} \mbox{ for all large negative $x$ and some $\delta_*>0$}.
\]
Let $\bar M:=M+3/\delta_*$; then at $r=c_*t-\bar M\log t$ we have
$$
\begin{aligned}
u(r,t)\ge& (1-Me^{-\delta t})\Phi_{c_*}(-3\delta_*^{-1}\log t)\ge (1-Me^{-\delta t})(1-e^{-3\log t})
\\
=& (1-Me^{-\delta t})(1-t^{-3})\ge 1-t^{-2}
\end{aligned}
$$
for all large $t$, say $t\geq \bar T\geq T$.
\end{proof}

\begin{Lemma}
\label{lem:l-bd}
There exist positive constants $T$ and $C$ such that
\[\mbox{
$h(t)\ge c_*t-(N-1)c^*\log t-C$ for $t\geq T$.}
\]
\end{Lemma}
\begin{proof}
In this case we perform comparison in the region $r\ge c_*t-\bar M\log t$, $t\geq T$, where $\bar M$ is as in the previous lemma, and $T$ is to be determined. The subsolution that we will use is given by
$$
\begin{aligned}
&w(r,t)=g(t)\Phi\Big(r-\underline{h}(t);\alpha\big(c_*-\frac{(N-1)c^*}t\big)\Big),\\
&g(t)=1-\frac{\log (t-t_0)}{(t-t_0)^2},\qquad \underline{h}(t)=c_*t-(N-1)c^*\log t+\frac{b\log (t-t_0)}{t-t_0}-C,
\end{aligned}
$$
where the function $\alpha(c)$ is the inverse of $c(\alpha)$, which is a well defined $C^2$ function (due to $c'(\alpha)\not=0$ and $c(\cdot)\in C^2(\mathbb R)$) with $\alpha(c_*)=0$, and the positive constants $t_0$, $b$ and $C$ will be chosen later.

It is trivially checked that $w(r,t)\le u(r,t)$ for $r=c_*t-\bar M\log t$ and $t\ge \max\{\bar T, t_0+3\}$, with $\bar T$ given in the previous lemma. Let us note that the functions $({\log t})/{t^2}$ and $(\log t)/t$ are decreasing for $t\geq 3$.

For $r\in(c_*t-\bar M\log t,\underline{h}(t))$,
$$
\begin{aligned}
\mathcal{L}w=&g'(t)\Phi-\underline{h}'(t)g(t)\Phi'
+\frac{(N-1)c^*\alpha'(c_*-\frac{(N-1)c^*}{t})g(t)}{t^2}\Phi_\alpha
\\
&-g^m(t)(\Phi^m)''-\frac{(N-1)g^m(t)}{r}(\Phi^m)'-g\Phi(1-g\Phi).
\end{aligned}
$$
Here we have omitted the obvious argument in $\Phi$, $\Phi^m$ and their derivatives for simplicity, and have denoted the derivative of $\Phi(x;\alpha)$ with respect to $\alpha$ by $\Phi_\alpha$ and the derivatives with respect to $x$ by primes.

Substituting
$$
(\Phi^m)''=-\Big(c_*-\frac{(N-1)c^*}{t}\Big)\Phi'-\alpha\Big(c_*-\frac{(N-1)c^*}{t}\Big)(\Phi^m)'-\Phi(1-\Phi)
$$
into the above expression, we obtain
$$
\begin{aligned}
\mathcal{L}w=&
\Phi_\alpha\underbrace{\frac{(N-1)c^*\alpha'(c_*-\frac{(N-1)c^*}{t})g(t)}{t^2}}_{\mathcal{A}}
\\
&-g\Phi'\left[\Big(\underbrace{(1-g^{m-1}(t))\big(c_*-\frac{(N-1)c^*}{t}\big)+\frac{b}{(t-t_0)^2}-\frac{b\log (t-t_0)}{(t-t_0)^2}}_{\mathcal{B}}\Big)\right.
\\
&\hspace{2cm} \left.+\underbrace{mg^{m-1}(t)\Phi^{m-1}\Big(\frac{N-1}{r}-\alpha\big(c_*-\frac{(N-1)c^*}{t}\big)}_{\mathcal{C}}\Big)\right]\\
&+\Phi\Big[\underbrace{\frac{2\log (t-t_0)}{(t-t_0)^3}-\frac1{(t-t_0)^3}+g^m(t)(1-\Phi)-g(t)(1-g(t)\Phi)}_{\mathcal{D}}\Big].
\end{aligned}
$$
Since $\alpha'(c_*-\frac{(N-1)c^*}{t})=\alpha'(c_*)+o(1)$, we clearly have
\[
\mathcal{A}=O(1/t^2).
\]
 For large $t$ and $t-t_0$, say $t\geq \tilde T$ and $t-t_0\geq \tilde M_0\geq 3$,
 and $b\geq  8(m-1)c_*$, we easily see that
$$
\mathcal{B}\le \frac{2(m-1)c_*\log (t-t_0)}{(t-t_0)^2}-\frac{b\log (t-t_0)}{2(t-t_0)^2}<-\frac{b\log t}{4t^2}.
$$
As for $\mathcal{C}$, we observe that
$$
\alpha\Big(c_*-\frac{(N-1)c^*}{t}\Big)=-\frac{\alpha'(c_*)(N-1)c^*}{t}+O(1/t^2)=\frac{N-1}{c_* t}+O(1/t^2)
$$
since $\alpha'(c_*)=-1/(c^*c_*)$. Hence, in view of $r\geq c_*t-\bar M \log t$ and $g, \Phi<1$, we have for large $t$,
$$
\begin{aligned}
\mathcal{C}\leq&\; m g^{m-1}(t) \Phi^{m-1}\left(\frac{N-1}{c_*t-\bar M\log t}-\alpha\big(c_*-\frac{(N-1)c^*}{t}\big)\right)\\
=&\;m g^{m-1}(t) \Phi^{m-1}\left(\frac{N-1}{c_*t}\big[1+\frac{\bar M\log t}{c_*t}+O(\log^2 t/t^2)\big]-\frac{N-1}{c_*t}+O(1/t^2)\right)
\\
=&\;m g^{m-1}(t) \Phi^{m-1}\left(\frac{(N-1)\bar M\log t}{(c_*t)^2}+O(1/t^2)\right)\\
\leq&\; \frac{2m(N-1)\bar M\log t}{(c_*t)^2}.
\end{aligned}
$$
As for the last term $\mathcal D$, using
\[
g^m(1-\Phi)-g(1-g\Phi)=-g\big[(1-g^{m-1})(1-\Phi)+(1-g)\Phi\big],
\]
 we get that
$$
\mathcal{D}\le \frac{2\log (t-t_0)}{(t-t_0)^3}-\frac1{(t-t_0)^3}-\frac{\mu_m\log (t-t_0)}{(t-t_0)^2}\le
-\frac{\mu_m\log (t-t_0)}{2(t-t_0)^2}\le -\frac{\mu_m\log t}{2t^2},
$$
provided that $t-t_0\geq M_0$ for some constant $M_0\geq \tilde M_0$ depending on $\mu_m$ given by
\[
\mu_m:=\min\Big\{\frac{m-1}8,\frac14\Big\}>0.
\]

We thus obtain, for large $t$ and  $r\in (c_*t-\bar M\log t,\;\underline h(t))$, with $0<t_0\leq t-M_0$,
\[
\mathcal L w\leq \Phi_\alpha \cdot O(1/t^2)-  g(t)\Phi' \left[-\frac{b\log t}{4t^2}+\frac{2m(N-1)\bar M\log t}{(c_*t)^2} \right]
-\Phi\cdot\frac{\mu_m\log t}{2t^2}.
\]
Using once more that $0<\Phi<1$, $1/2<g(t)<1$, we conclude, taking $b\geq \frac{16m(N-1)\bar M}{c_*^2}$, that
\[
\mathcal{L}w\le \Phi_\alpha \cdot O(1/{t^2})+\left(\frac{b}{16}\Phi'-\frac{\mu_m}2 \Phi\right)\frac{\log t}{t^2}
\]
and
$$
(m-1)\Phi^{m-2}\mathcal{L}w\le (\Phi^{m-1})_\alpha O\big(\frac 1{t^2}\big)+\left(\frac{b}{16}(\Phi^{m-1})'-(m-1)\frac{\mu_m}2 \Phi^{m-1}\right)\frac{\log t}{t^2}.
$$
Thus, since
\[\frac{b}{16}(\Phi^{m-1})'-(m-1)\frac{\mu_m}2 \Phi^{m-1}\leq - \nu \mbox{ in $\mathbb{R}_-$ for some $\nu>0$,}
\]
 and $(\Phi^{m-1})_\alpha=O(1)$ by Lemma~\ref{lem:phi.alpha.order.1}, we finally get
\[
\mbox{$\mathcal{L}w\leq 0$ for all large $t$, say $t\geq T\geq \tilde T$,
and $r\in (c_*t-\bar M\log t,\;\underline h(t))$,}
\]
 provided that we choose
\[
b\geq \max\left\{8(m-1)c_*,   \frac{16m(N-1)\bar M}{c_*^2}   \right\} \mbox{ and } t_0\in (0, T-M_0].
\]

With $b$, $T$ and $t_0$ chosen as above,
for $r\geq c_*T-\bar M\log T$, we have
\[
r-\underline h(T)\geq (c_*-\bar M)\log T-\frac {b\log (T-t_0)}{T-t_0}+C=0
\]
if we take
\[C:=-(c_*-\bar M)\log T+\frac {b\log (T-t_0)}{T-t_0}.
\]
It follows that
\[
w(r,T)=0\leq u(r, T) \mbox{ for } r\geq c_*T-\bar M\log T.
\]

We may now apply the comparison principle to conclude that $h(t)\geq \underline{h}(t)$, which yields immediately the desired estimate.
\end{proof}

%

\begin{Lemma}\label{lem:up-bd}
$h(t)\le c_* t-(N-1)c^*\log t+C$.
\end{Lemma}
\begin{proof}
In this case we compare in the set $r\ge0$, using as supersolution
$$
\begin{aligned}
&w(x,t)=g(t)\Phi\Big(r-\overline{h}(t);\alpha\big(c_*-\frac{(N-1)c^*}t\big)\Big),\\
&g(t)=1+\frac{\log (t-t_0)}{(t-t_0)^2},\qquad \overline{h}(t)=c_*t-c^*\log t-\frac{b\log (t-t_0)}{t-t_0}+C.
\end{aligned}
$$
For all large $t$,
$$
w(0,t)\ge \big(1+\frac{\log t}{t^2}\big)(1-e^{-\delta_1 t})\ge\big(1+\frac{\log t}{2t^2}\big)\ge 1+Me^{-\delta t}\ge u(0,t).
$$
We note that $w\equiv 0$ for $r>\overline h(t)$ and hence $\mathcal L w=0$ holds trivially in this range. We next show that
$\mathcal L w\ge 0$ in the range $r\in [0, \overline h(t)]$ for all large $t$.
Indeed, for $r\le \overline{h}(t)$, a computation very similar to that in Lemma \ref{lem:l-bd} yields $\mathcal{L}w\ge 0$ if $t$ is large enough, say $t\geq T$, and $t_0\in (0, T-M_0]$ for some $M_0\geq 3$. The initial data are ordered easily by taking $C$ large. We conclude that 
for $t\geq T$, $u\leq w$ and hence $h(t)\le \overline{h}(t)$, from which the result follows.
\end{proof}

Clearly Theorem \ref{thm:shift} follows directly from Lemmas \ref{lem:l-bd} and \ref{lem:up-bd}. Let us note that, from their proofs, we also have, for all large $t>0$ and some $C>0$,
\begin{equation}
\label{eq:u-bd}
\left\{\begin{array}{ll}
u(r,t)\geq \Phi^-(r-c_*t+(N-1)c^*\log t+C,t)\quad\text{for } r\ge c_*t-\bar M\log t,\\[8pt]
u(r,t)\leq \Phi^+(r-c_*t+(N-1)c^*\log t-C,t)\quad\text{for } r\ge0,
\end{array}\right.
\end{equation}
where
\begin{equation}
\label{eq:barriers}
\Phi^\pm(r,t)=\big(1\pm\frac{\log t}{t^2}\big)\Phi\big(r;\alpha(c_*-\frac{(N-1)c^*}{t})\big).
\end{equation}
We note that
\[
\lim_{t\to\infty}  \Phi^\pm(r,t)=\Phi_{c^*}(r).
\]
Hence, roughly speaking, for large times $u$ is trapped betweem two traveling wave solutions.

These estimates will play important roles in our proof of convergence in Section 6. Moreover, we note that Corollary \ref{cor:1}  follows directly from
\eqref{eq:u-bd} and Theorem \ref{thm:shift}.

\section{A uniform bound for the flux}
\label{sect:basic properties}
\setcounter{equation}{0}

The main aim of this section is to obtain a uniform bound for the flux $(u^m)_r$ for large times.  Here we always assume that
$T_0>0$ is chosen so that for every $t\geq T_0$,
\[
u(r, t)>0 \mbox{ for } r\in [0, h(t)),\; u(r,t)=0 \mbox{ for } r\geq h(t).
\]
 We will denote
\[
v(r,t):=\frac{m}{m-1}u^{m-1}(r,t).
\]
Then clearly
\[
v_t=(m-1)v \big(v_{rr}+\frac{N-1}{r}v_r\big)+v_r^2+f(v)\quad\text{for } r\in (0,h(t)),
\]
with
\[ f(v)=(m-1)v\left[1-\big(\frac{m-1}{m}v\big)^{\frac{1}{m-1}}\right].
\]


The following estimate, which was proved in~\cite{Perthame-Quiros-Vazquez-2014}, will be  frequently used in our analysis to follow:
\begin{equation}
\label{lem:pqv}
w(r,t):=\Delta v+f(v)\geq W(t):=\frac{-k \text{\rm e}^{-(m-1)kt}}{1-\text{\rm e}^{-(m-1)kt}},
\end{equation}
where
\[ k=k_F:=\min_{\xi\in [0,\frac{m}{m-1}]}[F(\xi)-\xi F'(\xi)]=\min\big\{1, \frac{1}{m-1}\big\}>0,\; F(\xi):=\frac{f(\xi)}{(m-1)\xi}.
\]
 Though the estimate holds only in the sense of distributions in $Q$, it holds in a pointwise sense in the set $\{(r,t): r<h(t),\,t>0\}$, where $u$ is positive and smooth.

\begin{Lemma}
\label{lem:v_r:h(t)}
$\lim_{r\to h(t)^-}v_r(r,t)=\sigma\in (-\infty, 0]$ exists for every $t>T_0$.
\end{Lemma}
\begin{proof} We fix $t>T_0$ and complete the proof in two steps.

\noindent\textsc{Step 1. } $\limsup_{r\to h(t)^-}v_r(r,t)<+\infty$.

Arguing indirectly we assume that
\[
\lim_{n\to \infty}v_r(a_n, t)=+\infty
\]
 along some sequence $\{a_n\}$ increasing to $h(t)$
as $n\to\infty$. Since $v(r,t)>0$ for $r\in [0, h(t))$ and $v(h(t), t)=0$, clearly
$\liminf_{r\to h(t)^-}v_r(r,t)\leq 0$. Therefore by passing to a subsequence of $\{a_n\}$ if necessary,
we can find another sequence $\{b_n\}$ satisfying
\[
\lim_{n\to\infty} v_r(b_n, t)=0,\; v_r(r,t)\geq 0 \mbox{ for } r\in (a_n, b_n),\; a_n<b_n<a_{n+1} \mbox{ for } n\geq 1.
\]
By estimate~\eqref{lem:pqv} and the uniform boundedness of $f(v)$, for each $t_0>0$ there exists $K=K_{t_0}$ such that
\begin{equation}
\label{4.2}
v_{rr}(r,t)+\frac{N-1}{r}v_r(r,t)\geq -K \mbox{ for } r\in (0, h(t)),\; t\geq t_0.
\end{equation}
We thus obtain
\[
\int_{a_n}^{b_n}\left[v_{rr}(r,t)+\frac{N-1}{r}v_r(r,t)\right]\,dr\geq -K(b_n-a_n)\to 0 \mbox{ as } n\to\infty.
\]
On the other hand,
\[
\int_{a_n}^{b_n}\left[v_{rr}(r,t)+\frac{N-1}{r}v_r(r,t)\right]\,dr\leq v_r(b_n, t)-v_r(a_n, t)+\frac{N-1}{a_n}\big[ v(b_n,t)-v(a_n, t)\big]\to-\infty
\]
as $n\to\infty$, since $\|v(\cdot, t)\|_\infty\leq \|v_0\|_\infty+\frac{m}{m-1}$. This contradiction completes Step 1.

\noindent\textsc{Step 2. } $\lim_{r\to h(t)^-}v_r(r,t)\in (-\infty, 0]$ exists.

Since $r^{-1}v_r(r, t)$ can be regarded as a smooth function of $r$ in $[0, h(t))$, the conclusion proved in Step 1 implies the existence of some
$M>0$ such that $\frac{N-1}{r}v_r(r,t)\leq M$ for $r\in [0, h(t))$, $t\geq T_0$. Therefore, in view of \eqref{4.2} there exists $K=K_{T_0}$ such that
\[ v_{rr}(r,t)\ge -K-\frac{N-1}{r}v_r\geq -K-M \mbox{ for } r\in [0, h(t)),\; t\geq T_0.
\]
This implies that $v_r(r,t)+(K+M)r$ is increasing in $r$ and hence $\lim_{t\to h(t)^-}v_r(r,t)>-\infty$ exists.
We already observed in the proof of Step 1 that $\liminf_{t\to h(t)^-}v_r(r,t)\leq 0$. It follows that
$\lim_{r\to h(t)^-}v_r(r,t)=\sigma\in (-\infty, 0]$.
\end{proof}

\begin{Theorem}
\label{thm:v_r-bd}
There exist $T_1\geq T_0$ and $C>0$ such that
$|(u^m)_r|,\, |v_r|\leq C $ for $t>T_1$ and $r\in [0, h(t))$.
\end{Theorem}
\begin{proof} 
It suffices to show the estimate for $|v_r|$, as $|(u^m)_r|=u|v_r|$ and $u\geq 0$ is uniformly bounded.
By Theorem \ref{thm:shift} and \eqref{eq:u-bd}, there exist $T'_0\geq T_0$ and $L>0$ such that
\[
v(r,t)\geq 1/2 \mbox{ for } t\geq T'_0,\; r\in [0, h(t)-L].
\]
For clarity we divide the arguments below into several steps.
\smallskip

\noindent\textsc{Step 1. } There exists $M_1>0$ such that $|v_r(r,t)|\leq M_1$ for $t\geq T_1:=T'_0+1$ and $r\in [0, h(t)-L-1]$.

Indeed, for $(r,t)$ in this range, since $1/2\leq v(r,t)\leq \max\{\|v_0\|_\infty, \frac m{m-1}\}$, by standard interior parabolic estimates we obtain $\|v\|_{C^{2,1}([R, R+1]\times [T, T+1])}\leq C$ for some $C$ independent of $R$ and $T$ as long as
$[R, R+1]\subset [0, h(t)-L-1]$ and $[T, T+1]\subset [T_0'+1, \infty)$. The required estimate follows directly from this.
\smallskip

\noindent\textsc{Step 2. } There exists $M_2>0$ such that $|v_r(h(t)^-, t)|\leq M_2$ for $t\geq T_1$.

Since $v_r(h(t)^-, t)\leq 0$, if the required estimate does not hold, then we can find a sequence $t_n\geq T_1$ such that
$\lim_{n\to\infty}v_r(h(t_n)^-, t_n)=-\infty$. Hence there exists $K_n\to+\infty$ and $0>r_n-h(t_n)\to 0$ such that
$
v_r(r_n, t_n)=-K_n$. We may assume that $K_n>1$ for all $n\geq 1$. Since $v_r(0, t_n)=0$ we can find a sequence $x_n\in (0, r_n)$ such that
\[
v_r(x_n, t_n)=-1, \; v_r(r,t_n)<-1 \mbox{ for } r\in (x_n, r_n).
\]
By  \eqref{4.2} there exists $M=M_{T_1}>0$ such that
\[
v_{rr}(r,t)+\frac{N-1}{r}v_r(r,t)\geq -M \mbox{ for } t\geq T_1,\; r\in [0, h(t)).
\]
It follows that
\[
v_{rr}(r,t_n)\geq -M-\frac{N-1}{r} v_r(r, t_n)\geq -M \mbox{ for } r\in (x_n, r_n).
\]
We thus obtain
\[
-K_n+1=v_r(r_n, t_n)-v_r(x_n, t_n)=\int_{x_n}^{r_n}v_{rr}(r,t_n)dr\geq -M(r_n-x_n).
\]
This implies $r_n-x_n\to\infty$ as $n\to\infty$ and hence
\[
v(x_n, t_n)=v(r_n, t_n)-\int_{x_n}^{r_n}v_r(r,t_n)dr\geq r_n-x_n\to\infty,
\]
a contradiction to the uniform boundedness of $v$. This completes the proof of Step 2.
\smallskip

\noindent\textsc{Step 3. } With $T_1$ as given in Step 1, there exists $M_3>0$ such that $|v_r(r,t)|\leq M_3$ for $t\geq T_1$ and $r\in [0, h(t))$.

Otherwise, in view of Steps 1 and 2,  there exists a sequence $\{(r_n, t_n)\}$ satisfying
\[
r_n\in (h(t_n)-L, h(t_n)),\; t_n\geq T_1 \mbox{ and } \lim_{n\to\infty}K_n:=|v_r(r_n,t_n)|= \infty.
\]
By passing to a subsequence we have two possibilities:
\[
{\rm (a) } \lim_{n\to\infty} v_r(r_n,t_n)=+\infty,\;\;\; {\rm (b)}  \lim_{n\to\infty} v_r(r_n,t_n)=-\infty.
\]
In case (a) we can find $y_n\in (r_n, h(t_n))$ such that for all large $n$, say $n\geq n_0$,
\[
v_r(y_n, t_n)=M_2+1,\; v_r(r, t_n)> M_2+1 \mbox{ for } r\in [r_n, y_n).
\]
Using \eqref{4.2} we obtain for $n\ge n_0$,
\[
\int^{y_n}_{r_n}\left[v_{rr}(r,t_n)+\frac{N-1}{r}v_r(r,t_n)\right]\,dr\geq -K(y_n-r_n)\ge -KL.
\]
On the other hand,
\[
\int^{y_n}_{r_n}\left[v_{rr}(r,t_n)+\frac{N-1}{r}v_r(r,t_n)\right]\,dr\leq v_r(y_n, t_n)-v_r(r_n, t_n)+\frac{N-1}{r_n}\big[ v(y_n,t)-v(r_n, t)\big]\to-\infty
\]
as $n\to\infty$, since $\|v(\cdot, t_n)\|_\infty\leq \|v_0\|_\infty+\frac{m}{m-1}$ and $r_n\geq h(T_1)-L>0$. 
This indicates that case (a) leads to a contradiction.

If case (b) happens, we can find $x_n\in (h(t_n)-L, r_n)$ such that for all large $n$, say $n\geq n_0$,
\[
v_r(x_n, t_n)=-(M_1+1),\; v_r(r, t_n)<-( M_1+1) \mbox{ for } r\in (x_n, r_n].
\]
Hence we obtain from \eqref{4.2} that
\[
v_{rr}(r,t_n)\geq -K \mbox{ for } r\in (x_n, r_n).
 \]
 It follows that
 \[
v_r(r_n, t_n)=v_r(x_n, t_n)+\int_{x_n}^{r_n}v_{rr}(r,t_n)dr\geq -(M_1+1) -K(r_n-x_n)\geq -(M_1+1)-KL,
\]
contradicting the fact that we are in case (b).
This completes the proof.
\end{proof}

\begin{Remark}
Using the above estimates  on $v_r$, it is not  difficult to show,  by adapting techniques in~\cite{Aronson-1986} and \cite{Knerr} for the one dimensional porous medium equation (without a source term), first that $h(t)$ is a Lipschitz continuous function, and then that
$h'(t)$ exists and $h'(t)=-v_r(h(t), t)$ for all large enough times.
\end{Remark}

\section{Identification of eternal solutions}
\label{sect:eternal}
\setcounter{equation}{0}

The following result, which identifies certain eternal solutions as wavefronts, will play a key role in the proof of convergence of $u(r,t)$  
along a subsequence of times. Such an identification result is well-known for the classical case $m=1$, but our current case $m>1$ is more difficult to treat due to the degeneracy of the differential operator. The proof here follows the strategy of the proof of a corresponding result in \cite{Du-Matsuzawa-Zhou-2015} for a free boundary problem, but the techniques here are completely different from  \cite{Du-Matsuzawa-Zhou-2015} as well as those for the $m=1$ case (see, e.g., \cite{BH-2007, Chen-1997, OM-1999}).

\begin{Theorem}
\label{thm:eternal}
Let $U(r,t)$ be a nonnegative weak solution to
\begin{equation}
\label{eq:eternal}
U_t=(U^m)_{rr}+c_*U_r+U(1-U)\quad\text{in }\mathbb{R}^2
\end{equation}
such that $U(r,t)$ is nonincreasing in $r$ and for some $C> 0$,
$$
\Phi_{c_*}(r+C)\leq U(r,t)\leq \Phi_{c_*}(r-C) \;\;\mbox{ for all } r\in\mathbb R.
$$
 Then there exists a constant $r_*\in[-C,C]$ such that $U(r,t)\equiv \Phi_{c_*}(r-r_*)$.
\end{Theorem}
The strategy of the proof is as follows. Let
\begin{equation}
\label{eq:def.R}
\begin{array}{l}
R^*=\inf\big\{R: U(r,t)\leq \Phi_{c_*}(r-R) \mbox{ for all } (r,t)\in \mathbb{R}^2 \big\},\\[8pt]
R_*=\sup\big\{R: U(t,r)\geq \Phi_{c_*}(r-R) \mbox{ for all } (t,r)\in\mathbb{R}^2 \big\}.
\end{array}
\end{equation}
Then $-C\le R_*\le R^*\le C$ and
\begin{equation}\label{5.3}
\Phi_{c_*}(r-R_*) \le U(r,t)\leq \Phi_{c_*}(r-R^*)\quad \mbox{for all } (r,t)\in\mathbb{R}^2.
\end{equation}
Hence the result will follow if we prove that $R_*=R^*$.

Since $U(\cdot,t)$ is nonincreasing for each $t\in\mathbb{R}$, there exists a value $G(t)$, the position of the unique interface (also called the free boundary) of $U(\cdot,t)$,  such that $U(r,t)>0$ for $r<G(t)$ and $U(r,t)=0$ for $r\ge G(t)$.

%


The first step is to show that $U(r,t)$ approaches $\Phi_{c_*}(r-R^*)$ for any fixed $r$ along some time sequence. 

\begin{Lemma}
	\label{lemma:M_r.above} Let $U$ be as in Theorem~\ref{thm:eternal}, $G(t)$ the function giving the position of its 
	free boundary at time $t$,  and $R^*$ as in~\eqref{eq:def.R}.  Then for any $r\in\mathbb{R}$,
	\begin{equation}
	\label{eq:definition.Mr}
	M_r:=\inf_{t\in\mathbb{R}}\big[\Phi_{c_*}(r-R^*)-U(r,t)\big]=0.
	\end{equation}
\end{Lemma}

\begin{proof} The result is obvious for $r\ge R^*$. Suppose that there exists some $\bar r<R^*$ such $M_{\bar r}=2\delta>0$. Let $\varepsilon>0$ be such that $\Phi_{c_*}(\bar r-(R^*-\varepsilon))=\Phi_{c_*}(\bar r-R^*)-\delta$. We will show that $U(r,t)\le  \Phi_{c_*}(\bar r-(R^*-\varepsilon))$ for $(r,t)\in\mathbb{R}^2$, thus contradicting the definition of $R^*$. 
	
We first prove the bound  for $r\le\bar r$. We consider the auxiliary problem
\begin{equation}
\label{bar-V}
\begin{cases}
\overline U_t=(\overline
U^m)_{rr}+c_*\overline U_r+\overline U(1-\overline U) &\text{for }  r<\bar r,\; t>0,\\
\overline U(\bar r,t)=\Phi_{c_*}(\bar r-(R^*-\epsilon)),&\text{for } t>0,\\
\overline U(r,0)=1, &\text{for } r<\bar r.
\end{cases}
\end{equation}
It
follows from the comparison principle that
\[
1\geq \overline U(r,t)\geq \overline U(r,t+h)\geq \psi(r):=\Phi_{c_*}(r-(R^*-\epsilon)) \mbox{ for all }
r<\bar r,\; t>0,\; h>0.
\] 
Thus, since $\overline U$ is nonincreasing in $t$, it has a limit as $t\to\infty$, which satisfies
\[
U^*(r):=\lim_{t\to\infty} \overline U(r,t)\geq
\psi(r)\text{ for all }r<\bar r.
\]
We will prove now that in fact $U^*(r)\equiv \psi(r)$.
Indeed, $U^*$ is a nonincreasing function that satisfies
\begin{equation}
\label{V*} 
[(U^*)^m]_{rr}+c_* U^*_r+U^*(1- U^*)\mbox{ for } r< \bar r,\;
U^*(-\infty)=1,\; U^*(\bar r)=\psi(\bar r),
\end{equation}
from where the identity $U^*\equiv \psi$ follows easily by an ODE argument. 
%
%

We now look at $U(r,t)$, which satisfies \eqref{eq:eternal}, and for any $t\in \mathbb{R}$,
\[
U(r,t)< 1,\; U(\bar r,t)\leq \Phi_{c_*}(\bar r-R^*)-2\delta\leq
\psi(\bar r).
\]
Therefore we can use the comparison principle to deduce that
\[
U(r,s+t)\leq \overline U(r,t) \mbox{ for all } t>0, r<\bar r, s\in\mathbb{R},
\]
or, equivalently,
\[
U(r,t)\leq \overline U(r,t-s) \mbox{ for all } t>s, r<\bar r,
s\in\mathbb{R}.
\]
Letting $s\to-\infty$ we obtain
\begin{equation}
\label{r<R0} U(r,t)\leq U^*(r)=\psi(r) \mbox{ for all }
r<\bar r, t\in\mathbb{R}.
\end{equation}

To prove that the bound holds for $r\in [\bar r,R^*]$ we perform an analogous argument, now using the auxiliary problem 
\begin{equation*}
\begin{cases}
\overline U_t=(\overline
U^m)_{rr}+c_*\overline U_r+\overline U(1-\overline U)\quad &\text{for }  r\in (\bar r,R^*),\; t>0,\\
\overline U(\bar r,t)=\psi(\bar r),\quad \overline U(R^*,t)=0&\text{for } t>0,\\
\overline U(r,0)=1\qquad &\text{for } r\in(\bar r,R^*).
\end{cases}
\end{equation*}
The solution to this latter problem converges monotonically to $\psi(r)$ as $t\to+\infty$. Comparing $U(r,t)$ with $\overline U(r,t-s)$ for $t>s$ and then letting $s\to-\infty$ we obtain that $U(r,t)\le \psi(r)$ for all $r\in[\bar r, R^*]$ and $t\in\mathbb{R}$. 

Summarizing,  $U(r,t)\le \psi(r)$ for all $(r,t)\in\mathbb{R}^2$, which yields the desired contradiction. 
\end{proof}

A similar argument, now using auxiliary subsolutions with zero initial datum, gives the following analogous result, whose detailed proof is omitted.

\begin{Lemma}
	\label{lemma:M_r.below} Let $U$ be as in Theorem~\ref{thm:eternal}, $G$ the function giving the position of its interface at each time,  and $R_*$ as in~\eqref{eq:def.R}. Then for any $r\in\mathbb{R}$,
	\begin{equation}
	\label{eq:definition.Mr}
m_r:=\inf_{t\in\mathbb{R}}\big[U(r,t)-\Phi_{c_*}(r-R_*)\big]=0.
	\end{equation}
\end{Lemma}	

We now prove that $U(r,t)$ and the function $G(t)$ giving its free boundary  converge along a sequence of time  to $\Phi_{c_*}(r-R^*)$ and $R^*$, respectively. 

\begin{Lemma}
	\label{limit-VG}
	Under the hypotheses of Lemma~\ref{lemma:M_r.above}, there exists a sequence $\{s_n\}\subset \mathbb{R}$ such that
	\[
	G(s_n)\to R^*,\; U(r, t+s_n)\to \Phi_{c_*}(r-R^*)
	\quad\mbox{as } n\to\infty
	\]
	uniformly for $(r,t)$ in compact subsets of $\mathbb{R}^2$.
\end{Lemma}
\begin{proof} Fix $r_0<R^*$. We know from Lemma~\ref{lemma:M_r.above} that $M_{r_0}=0$. Therefore, there are two possibilities:
	\begin{itemize}
		\item[(i)] $\Phi_{c_*}(r_0-R^*)=\sup_{t\in\mathbb{R}}U(r_0,t)$ is achieved at some finite $t=s_0$,
		\item[(ii)] $\Phi_{c_*}(r_0-R^*)>U(r_0,t)$ for all $t\in\mathbb{R}$ and $U(r_0,s_n)\to \Phi_{c_*}(r_0-R^*)$ along some unbounded sequence $s_n$.
	\end{itemize}
	
	In case (i), since $U(r,t)\leq \Phi_{c_*}(r-R^*)$ for $(r,t)\in\mathbb{R}^2$, with
	$U(r_0,s_0)=\Phi_{c_*}(r_0-R^*)>0$, we can apply the strong maximum principle (and the monotonicity of $U$ in $r$) to conclude that
	$U(r,t)\equiv\Phi_{c_*}(r-R^*)$ for all $(r,t)\in\mathbb{R}^2$. Thus the conclusion of the lemma holds by taking
	$s_n\equiv s_0$.
	
	In case (ii), we consider the sequence
	\[
	U_n(r,t)=U(r,t+s_n).
	\]
	The regularity results in~\cite{diBenedetto-1983,Ziemer-1982} guarantee that, by passing to a subsequence,
	\[
	U_n\to\tilde U \mbox{ locally  uniformly in } \mathbb{R}^2,
	\]
	where $\tilde{U}$  satisfies equation \eqref{eq:eternal}, and $\tilde{U}(r,t)\leq \Phi_{c_*}(r-R^*)$ for all $(r,t)\in\mathbb{R}^2$.
	Moreover,
	\[
	\tilde U(r_0,0)=\Phi_{c_*}(r_0-R^*).
	\]
	Hence we are back to case (i) and thus  $\tilde U(r,t)\equiv \Phi_{c_*}(r-R^*)$ for $(r,t) \in \mathbb{R}^2$, and so
	\[
	U(r, t+s_n)\to \Phi_{c_*}(r-R^*) \mbox{ locally  uniformly in } \mathbb{R}^2.
	\]

We now consider $G(s_n)$. From $U(r,t)\leq \Phi_{c_*}(r-R^*)$ for $(r,t)\in\mathbb{R}^2$, we clearly have $G(t)\leq R^*$.
If $\liminf _{n\to\infty} G(s_n)=\tilde R<R^*$, say $G(s_{n_k})\to \tilde R$ along some subsequence $\{s_{n_k}\}$ of $\{s_n\}$,
then, as $k\to\infty$, $0=U(G(s_{n_k}), s_{n_k})\to \Phi_{c_*}(\tilde R-R^*)>0$, a contradiction. Therefore we must have
$\liminf _{n\to\infty} G(s_n)=R^*$, which implies $\lim_{n\to\infty} G(s_n)=R^*$ due to $G(t)\leq R^*$ for all $t$.
\end{proof}

Next we prove an analogue of Lemma \ref{limit-VG} with $R^*$ replaced by $R_*$. However, an extra difficulty arises in the proof, due to the possibility that the solution may ``degenerate'', namely $U(r,t_n)\to\Phi_{c_*}(r-R_*)$ along some time sequence $t_n$, but $\limsup_{n\to\infty} G(t_n)>R_*$.
 To cope with this complication we introduce an additional technique called a ``blocking method'' (see details below). 

\begin{Lemma}
	\label{limit-VG.below}
	Under the hypotheses of Lemma~\ref{lemma:M_r.below}, there exists a sequence $\{\tilde s_n\}\subset \mathbb{R}$ such that
	\[
	G(\tilde s_n)\to R_*,\; U(r, t+\tilde s_n)\to \Phi_{c_*}(r-R_*)
	\quad\mbox{as } n\to\infty
	\]
	uniformly for $(r, t)$ in compact subsets of $\mathbb{R}^2$.
\end{Lemma}
\begin{proof}
Since Lemma~\ref{lemma:M_r.below} gives $m_r=0$ for all $r\in\mathbb{R}$, we can argue as in the proof of Lemma~\ref{limit-VG} to obtain the existence of a sequence $\{t_n\}$ such that $U(r,t+t_n)\to \Phi_{c_*}(r-R_*)$ uniformly on compact subsets of $\mathbb{R}^2$. However, we are not able to show $G(t_n)\to R_*$ as  in the proof of Lemma~\ref{limit-VG}, since $\limsup_{n\to\infty} G(t_n)>R_*$ does not lead to a contradiction as before.

To show the existence of a sequence $\{\tilde s_n\}$ having the properties stated in the lemma,  we need some new techniques. Clearly, for each $k\in\mathbb{N}$, there exists $n_k\in\mathbb N$ such that
\begin{equation}
\label{U-Phi}
|U(r, t_n+t)-\Phi_{c_*}(r-R_*)|\leq k^{-1} \mbox{ for } r, t\in [-2k, 2k],\; n\geq n_k.
\end{equation}
We show that, for each $j\in\mathbb N$, there exists $k_j\in\mathbb N$ so that
\begin{equation}
\label{G-R_*}
 \mbox{$ 0\le \min_{t\in [0, k]}G(t_{n_k}+t)-R_*\leq j^{-1}$ for all $k\geq k_j$.}
 \end{equation}
 
 If \eqref{G-R_*} is proved, then for each $j\in \mathbb N$ there exists $s_j\in [0, k_j]$ such that
 \[
 0\leq G(t_{n_{k_j}}+s_j)-R_*\leq j^{-1}.
 \]
 If we denote $\tilde s_j:=t_{n_{k_j}}+s_j$, then clearly $G(\tilde s_j)\to R_*$ as $j\to\infty$. Moreover, by 
 \eqref{U-Phi} we have
 \[
 |U(r, \tilde s_j+t)-\Phi_{c_*}(r-R_*)|\leq \frac 1{k_j} \mbox{ for } r,t\in [-k_j, k_j].
 \]
 Hence
 \[
 U(r, \tilde s_j+t)\to \Phi_{c_*}(r-R_*) \mbox{ as } j\to \infty \mbox{ uniformly for $(r,t)$ in compact subsets of } \mathbb R^2.
 \]
 Thus the conclusions of the lemma hold provided that \eqref{G-R_*} holds.
 
 So to complete the proof of the lemma it suffices to prove \eqref{G-R_*}. Suppose \eqref{G-R_*} does not hold; then
  there exists some $j\in\mathbb N$ and a sequence $k^i\to\infty$ as $i\to\infty$ such that
\begin{equation}
\label{eq:assumption.boundary}
G(t_{n_{k^i}}+t)-R_*\ge \delta:=j^{-1}>0 \mbox{ for all } t\in [0, k^i],\; i=1,2,...
\end{equation}
We now use \eqref{eq:assumption.boundary} to derive a contradiction by means of a ``blocking method''. 

We define
\[
u(r,t):=U(r-c_*t, t),\quad g(t):=c_*t +G(t).
\]
Then clearly $r=g(t)$ is the free boundary of $u$ and
\[
u_t=(u^m)_{rr}+u(1-u) \mbox{ in } \mathbb R^2.
\]
By \eqref{U-Phi} we have, for all large $k\in\mathbb N$,
\[
|u(r+c_*(t_{n_k}+t), t_{n_k}+t)-\Phi_{c_*}(r-R_*)|<\frac 1k \mbox{ for } (r,t)\in [R_*,R^*]\times [0, 2k].
\]
Thus, from~\eqref{eq:assumption.boundary} and the definition of $g$  it follows that, for all large positive integer $i$,
\[
\mbox{$[G(t_{n_{k^i}}+t)-\delta,G(t_{n_{k^i}}+t)]\subset [R_*, G(t_{n_{k^i}}+t)]$ for $t\in [0, k^i]$,}
\]
and
\begin{equation}
\label{u-small}
0<u(r+c_*(t_{n_{k^i}}+t), t_{n_{k^i}}+t)<\frac 1{k^i} \mbox{ for } r\in [R_*, G(t_{n_{k^i}}+t)),\; t\in [0, k^i].
\end{equation}
This will allow us to construct a supersolution for $(r,t)$ in a suitable range to show that, for all large $i$, the free boundary of $u$ has to grow very slowly when $t\in [t_{n_{k^i}}, t_{n_{k^i}}+ k^i]$, because it is blocked by the free boundary of the constructed supersolution, which by construction 
grows very slowly.
This slow growth of the free boundary of $u$  would easily induce a contradiction. We call this method a ``blocking method'', and will use it again in the next section.

The supersolution is given by
\[\bar u(r,t):=\left(\frac{m-1}{m}\right)^{\frac1{m-1}}
\text{\rm e}^{t}\,\Big[\epsilon^2\frac{\text{\rm e}^{(m-1)t}-1}{m-1}-\epsilon r\Big]_+^{\frac1{m-1}}
\]
with $\epsilon>0$. A straightforward computation shows that it
satisfies
\[
\bar u_t-(\bar u^m)_{rr}=\bar u \mbox{ for } (r,t)\in\mathbb{R}\times\mathbb{R}_+,
\]
with free boundary given by
\[
r=\bar h(t):=\epsilon\frac{\text{\rm e}^{(m-1)t}-1}{m-1}.
\]
Note that $\bar u$ is monotone both in space and time.  Set
$\delta^*=\frac 13\min\{\delta, c_*\}$,
 and then fix $\epsilon>0$ small so that the free boundary of the supersolution does not advance much in the time interval $[0,1]$, namely,
\[
\bar{h}(1)-\bar{h}(0)=\bar h(1)=\epsilon \frac{\text{\rm e}^{m-1}-1}{m-1}<\delta^*;
\]
and $\bar u(r,t)$ is not too big  at $r=-2\delta_*$ in this time interval:
\[
0<\bar u(-2\delta^*, 0)\leq \bar u(-2\delta^*, t)\leq \bar u(-2\delta^*, 1)<1 \mbox{ for } t\in [0,1].
\]
By
\eqref{u-small}, there exists $i_0>0$ large so that \begin{equation}\label{*}
0\leq u(r,t_{n_{k^i}}+t)<\bar u(-\delta^*,0) \mbox{ for } r\in [g(t_{n_{k^i}}+t)-\delta, g(t_{n_{k^i}}+t)],\; t\in [0, k^i],\; i\geq i_0.
\end{equation}
Fix $i\geq i_0$ and $t_0\in [0, k^i-1]$ and define
\[
\hat u(r, t):=\bar u(r-g(t_{n_{k^i}}+t_0)-\delta^*, t-t_0),\qquad \tilde u(r,t):=u(r, t_{n_{k^i}}+t),
\]
whose free boundaries are given respectively by
\[
r=\hat h(t):=\bar h(t-t_0)+g(t_{n_{k^i}}+t_0)+\delta^*,\qquad r=\tilde h(t):=g(t_{n_{k^i}}+t).
\]
Clearly, 
\[
\left\{
\begin{array}{l}
\displaystyle\hat u_t-(\hat u^m)_{rr}=\hat u,\\[8pt]
\displaystyle \tilde u_t-(\tilde u^m)_{rr}=\tilde u(1-\tilde u)\leq \tilde u,
\end{array}
\right.
\qquad\mbox{for } r\in \mathbb{R}_+,\; t\geq t_0,
\]
 Moreover,
\[
 \hat h(t_0)=\tilde h(t_0)+\delta^*,\qquad \hat h(t_0+1)=\bar h(1)+g(t_{n_{k^i}}+t_0)+\delta^*<\tilde h(t_0)+2\delta^*.
\]

We claim that
\[
\tilde h(t)\leq \hat h(t) \mbox{ for } t\in [t_0, t_0+1].
\]
If this claim is proved, then we have
\begin{equation}\label{ht0}
\tilde h(t_0+1)-\tilde h(t_0)\le\hat h(t_0+1)-\tilde h(t_0)<2\delta^*\leq\frac{2}{3}c_*.
\end{equation}
Because \eqref{ht0} holds for all $t_0\in [0, k^i-1]$, it follows that
\[
g(t_{n_{k^i}}+k^i)-g(t_{n_{k^i}})=\tilde h(k^i)-\tilde h(0)\leq \frac{2}{3}c_*k^i\;\; \mbox{ for } i\geq i_0.
\]
Therefore
\[
c_*k^i+G(t_{n_{k^i}}+k^i)-G(t_{n_{k^i}})\leq \frac 23 c_* k^i \mbox{ for } i\geq i_0,
\]
 which is an obvious contradiction
to $R_*\leq G(t)\leq R^*$ in $\mathbb R$ when $i$ is sufficiently large. So, to conclude the proof of \eqref{G-R_*} it suffices to prove the above claim.

Arguing indirectly we assume that the claim is not true. Then,
since $\tilde h(t_0)<\hat h(t_0)$, we can find some $t^0\in (t_0, t_0+1)$ so that $\tilde h(t^0)\in \big(\hat h(t^0), \hat h(t_0+1)\big)$.
Therefore, for $t\in [t_0, t^0]$ we have 
\[
\tilde h(t)-\tilde h(t_0)+\delta^*<\hat h(t_0+1)-\tilde h(t_0)+\delta^*<3\delta^*\leq \delta,
\]
 and so by \eqref{*}, for such $t$,
\begin{equation}\label{**}
\tilde u(\tilde h(t_0)-\delta^*, t)<\bar u(-\delta^*, 0)<\bar u(-2\delta^*,0)
=\hat u(\tilde h(t_0)-\delta^*, t_0) \smallskip\leq \hat u(\tilde h(t_0)-\delta^*, t).
\end{equation}
Choose $R_0>\max\left\{\hat h(t_0+1), \tilde h(t_0+1)\right\}$ so that
\begin{equation}\label{***}
\tilde h(t), \hat h(t) <R_0 \mbox{ for } t\in [t_0, t^0],
\end{equation}
and then compare $\tilde u(r, t)$ with $\hat u(r,t)$ in the region $(r, t)\in [\tilde h(t_0)-\delta^*, R_0]\times [t_0, t^0]$.
On the parabolic boundary of this region, by \eqref{*}, \eqref{**} and \eqref{***}, we have $\tilde u\leq \hat u$, and hence we can use the comparison principle to conclude that
$\tilde u(r, t)\leq \hat u(r, t)$ in this region. In particular, $\tilde  u(\hat h(t^0), t^0)\leq \hat u(\hat h(t^0), t^0)=0$.
On the other hand, from $\tilde h(t^0)>\hat h(t^0)$ we obtain $\tilde u(\hat h(t^0), t^0)>0$. This contradiction proves our claim, and thus \eqref{G-R_*} is proved.
\end{proof}

We are now ready to show that the upper and the lower optimal barriers coincide.
\begin{Lemma}
\label{G=const}
Let $U$ be as in Theorem~\ref{thm:eternal}, and $R_*$ and $R^*$ as in~\eqref{eq:def.R}. Then
$R_*=R^*$.
\end{Lemma}
\begin{proof}
From  Lemma 5.4 we obtain
\[
\lim_{n\to\infty}U(r, s_n)=\Phi_{c_*}(r-R^*) \mbox{ in } C_{loc}(\mathbb R),\; \lim_{n\to\infty} G(s_n)=R^*.
\]
Since $U(r, s_n)=\Phi_{c_*}(r-R^*)\equiv 0$ for $r\geq R^*+C$,  and $\Phi_{c_*}(r-R^*)\to 1$, $U(r, s_n)\to 1$ as $r\to -\infty$ uniformly in $n$ (due to \eqref{5.3}),
it follows easily that,  for any small $\delta>0$, we can find $n=n_\delta$ sufficiently large so that
\[
(1-\delta)\Phi_{c_*}(r-R^*+\delta )\leq U(r, s_n)\leq (1+\delta)\Phi_{c_*}(r-R^*-\delta) \mbox{ for } r\in \mathbb R.
\]
We are now in a position to use Lemma 3 of \cite{Biro-2002} to conclude that,  for $t>0$ and $r\in \mathbb R$,
\[
(1-\delta)\Phi_{c_*}(r-R^*+\omega_1(\delta) )\leq U(r, s_n+t)\leq (1+\delta)\Phi_{c_*}(r-R^*-\omega_2(\delta)),
\]
with 
\[
0<\omega_1(\delta)\leq L\delta+\frac{\delta}{1-\delta},\; 0<\omega_2(\delta)\leq L\delta+\delta,
\]
where the constant $L>0$ is independent of $\delta$. It follows that
\[
R^*-\omega_1(\delta)\leq G(s_n+t)\leq R^*+\omega_2(\delta) \mbox{ for } t>0.
\]

Similarly, from  Lemma 5.5 we obtain
\[
\lim_{n\to\infty}U(r, \tilde s_n)=\Phi_{c_*}(r-R_*) \mbox{ in } C_{loc}(\mathbb R),\; \lim_{n\to\infty} G(\tilde s_n)=R_*.
\]
We may now use Lemma 3 of \cite{Biro-2002} to analogously deduce the existence of some $m=m_{\delta}$ large such that
\[
R_*-\omega_1(\delta)\leq G(\tilde s_m+t)\leq R_*+\omega_2(\delta) \mbox{ for } t>0.
\]
Thus for any $t_0>\max\big\{s_n, \tilde s_m\big\}$, we have
\[
R^*-\omega_1(\delta)\leq G(t_0)\leq R_*+\omega_2(\delta),
\]
which gives
\[
0\leq R^*-R_*\leq \omega_1(\delta)+\omega_2(\delta).
\]
Letting $\delta\to 0$ we immediately obtain
 $R_*=R^*$.
\end{proof}
As mentioned earlier, Theorem \ref{thm:eternal} follows immediately from $R_*=R^*$.

\section{Convergence and some remarks}
\label{sect:convergence}
\setcounter{equation}{0}

In this section, we prove the convergence results on $u$ and $h$ stated in Theorem \ref{thm:main}. The idea is to obtain first convergence along a sequence of time $t_n\to\infty$, and then refine the sub- and supersolutions in the proof of Lemmas 3.4 and 3.5 
to show that once $u(t,r)$ is close to $\Phi_{c_*}(r-c_*t+(N-1)c^*\log t-r_0)$ at some time $t_n$, then it remains close for all later time.
These are done in subsections 6.1 and 6.2. In subsection 6.3, we will comment on a gap in the proof of the main result (Theorem 4) in \cite{Biro-2002} and indicate how the gap can be fixed by techniques developed here; we will also give a version of Theorem 1.1 in dimension one but without requiring $u_0$ to be symmetric.

\subsection{Convergence along a time sequence}  In this subsection we prove the following result.

\begin{Proposition}\label{prop6.1}
There exist a sequence $\{t_k\}$ with $\lim_{k\to\infty} t_k=\infty$ and a constant $r_0\in\mathbb{R}$ such that
$$
\lim_{k\to\infty}u(r+c_*t_k-(N-1)c^*\log t_k)=\Phi_{c^*}(r-r_0)
$$
uniformly in $r\ge a$, for any $a\in\mathbb{R}$, and
$$
\lim_{n\to\infty} \big[h(t_k)-c_*t_k+(N-1) c^*\log t_k\big]=r_0.
$$
\end{Proposition}

To prove this proposition, we prepare several lemmas.
\begin{Lemma}
\label{lem:convergence.solution.subsequence}
Given any sequence of times $\{t_n\}$ such that
$\lim_{n\to\infty} t_n=\infty$, there exist a subsequence $\{t_{n_k}\}$ and a constant $r_0\in\mathbb{R}$ such that
$$
\lim_{k\to\infty}u(r+c_*t_{n_k}-(N-1)c^*\log t_{n_k},t_{n_k})=\Phi_{c^*}(r-r_0)
$$
uniformly in $r\ge a$, for any $a\in\mathbb{R}$.
\end{Lemma}
\begin{proof}
Let $k(t):=c_*t-(N-1)c^*\log t$.
By  Theorem \ref{thm:shift} and \eqref{eq:u-bd}, there exist $T>0$ such that
\begin{align*}
\tilde{u}(r,t):=u(r+k(t),t),\quad \tilde{h}(t):=h(t)-k(t)
\end{align*}
satisfy, for all $r\ge -[\bar M-(N-1)c^*]\log t$ and $t\geq T$,
\begin{equation}\label{C-3C}
-C\leq \tilde{h}(t)\leq C,\quad \Phi^-(r+C, t)\leq \tilde u(r,t)\leq \Phi^+(r-C, t),
\end{equation}
where the barriers $\Phi^\pm$ are as in~\eqref{eq:barriers}.
Moreover,
\[
\tilde{u}_t=(\tilde{u}^m)_{rr}+\frac{N-1}{r+k(t)}(\tilde{u}^m)_r+\frac{c_*-(N-1)c^*}{t}\tilde{u}_r+\tilde u(1-\tilde u) \mbox{ for } r>-k(t),\ t>T.
\]
Let $\{t_n\}$  be an arbitrary sequence such that $t_n\to\infty$. Define
\[
 \tilde{u}_n(r,t)=\tilde{u}(r,t+t_n),\; \tilde{h}_n(t)=\tilde{h}(t+t_n).
\]
Thanks to the regularity results of~\cite{diBenedetto-1983,Ziemer-1982}, we know that there exists a subsequence, still denoted by itself for convenience, such that as $n\to\infty$,
\[
\tilde u_n(r,t)\to U(r,t) \mbox{ in } C_{\rm loc}(\mathbb R^2),
\]
where $U$ is  nonincreasing in $r$ and satisfies 
$$
U_t=(U^m)_{rr}+c_*U_r+U(1-U)\quad\text{in }\mathbb{R}^2, \qquad \Phi_{c_*}(r+C)\leq \tilde U(r,t)\leq \Phi_{c_*}(r-C)
$$
in a weak sense.
Here we have used  the uniform bound for $(\tilde u_n^m)_r$ proved in Theorem~\ref{thm:v_r-bd}, and the monotonicity of $u(r, t)$ for $r\geq h_0$.
Using now the identification of eternal solutions from Section~\ref{sect:eternal}, we conclude that $U(r,t)=\Phi_{c_*}(r-r_0)$ for some $r_0\in[-C,C]$ and hence the wanted result, by taking $t=0$.
\end{proof}


The above convergence result implies in particular that $\liminf_{n\to\infty}\tilde{h}(t_n)\ge r_0$. Unfortunately it does not follow automatically that $\lim_{n\to\infty}\tilde{h}(t_n)= r_0$. To prove the validity of this identity,
 we will need to prove a certain non-degeneracy result for $u(r,t)$  close to its free boundary $r=h(t)$. This is done in the following two lemmas.
\begin{Lemma}\label{sigma(r)}
For $r>0$ we have
\[
\sigma(r):=\limsup_{t\to\infty}u(h(t)-r, t)>0.
\]
\end{Lemma}
\begin{proof}
We prove the conclusion by an indirect argument.  Suppose there exists $\delta>0$ such that $\sigma(\delta)=0$.
Then clearly $\lim_{t\to\infty} u(h(t)-\delta, t)=0$. Due to the monotonicity of $u(r,t)$ in $r$ for $r>h(0)$ we obtain
\begin{equation}\label{v-small}
\lim_{t\to\infty}\left[\max_{r\in [h(t)-\delta, h(t)]}u(r,t)\right]=0.
\end{equation}

We are now in a position to produce a contradiction by using the ``blocking method'' in the proof of Lemma \ref{limit-VG.below}, with some simple variations.
Define $\bar u(r,t)$, $\bar h(t)$ and $\delta^*$ as in the proof of Lemma \ref{limit-VG.below}.
 We then fix $\epsilon>0$ small, as before, so that
\[
\bar{h}(1)-\bar{h}(0)=\bar h(1)=\epsilon \frac{\text{\rm e}^{m-1}-1}{m-1}<\delta^*
\]
 and
\[
0<\bar u(-2\delta^*, 0)\leq \bar u(-2\delta^*, t)\leq \bar u(-2\delta^*, 1)<1 \mbox{ for } t\in [0,1].
\]
By
\eqref{v-small}, there exists $T>0$ large so that
\begin{equation}\label{*1}
0\leq u(r,t)<\bar u(-\delta^*,0) \mbox{ for } r\in [h(t)-\delta, h(t)],\; t\geq T.
\end{equation}
Fix $t_0\geq T$ and let
\[
\hat u(r, t)=\bar u(r-h(t_0)-\delta^*, t-t_0),\; \hat h(t)=\bar h(t-t_0)+h(t_0)+\delta^*.
\]
Clearly, since $(\hat u^m)_r\le 0$,
\[
\left\{
\begin{array}{l}
\displaystyle\hat u_t-(\hat u^m)_{rr}- \frac{N-1}r (\hat u^m)_r\ge\hat u,\\[8pt]
\displaystyle u_t-(u^m)_{rr}- \frac{N-1}r ( u^m)_r=u(1-u)\leq u,
\end{array}
\right.
\qquad\mbox{for } r\in \mathbb{R}_+,\; t\geq t_0,
\]
and
\[
\hat u(\hat h(t), t)=0,\; \hat h(t_0)=h(t_0)+\delta^*,\; \hat h(t_0+1)=\bar h(1)+h(t_0)+\delta^*<2\delta^*+h(t_0).
\]

We claim that
\[
h(t)\leq \hat h(t) \mbox{ for } t\in [t_0, t_0+1].
\]
If this claim is proved, then we have
\begin{equation}\label{ht01}
h(t_0+1)-h(t_0)<\hat h(t_0+1)-h(t_0)<2\delta^*\leq\frac{2}{3}c_*.
\end{equation}
Because \eqref{ht01} holds for all $t_0>T$, the average speed of $h(t)$  is thus no bigger than $\frac{2}{3}c_*$, which is an obvious contradiction
to \eqref{C-3C}.

 So to complete the proof, it suffices to prove the above claim, which can be done the same way as in he proof of Lemma \ref{limit-VG.below}.
 We give the details below for convenience of the reader. Arguing indirectly we assume that the claim is not true. Then
since $h(t_0)<h(t_0)+
\delta^*=\hat h(t_0)$, we can find some $t^0\in (t_0, t_0+1)$ so that
\[
        h(t^0)\in \big(\hat h(t^0), \hat h(t_0+1)\big).
\]
Therefore, for $t\in [t_0, t^0]$ we have $h(t)-h(t_0)+\delta^*<\hat h(t_0+1)-h(t_0)+\delta^*<3\delta^*\leq \delta$ and so by \eqref{*1}, for such $t$,
\begin{equation}\label{**1}
u(h(t_0)-\delta^*, t)<\bar u(-\delta^*, 0)<\bar u(-2\delta^*,0)=\hat u(h(t_0)-\delta^*, t_0)\leq \hat u(h(t_0)-\delta^*, t).
\end{equation}
Choose $R_0>\max\{\hat h(t_0+1), h(t_0+1)\}$ so that
\begin{equation}\label{***1}
h(t), \hat h(t) <R_0 \mbox{ for } t\in [t_0, t^0],
\end{equation}
and then compare $u(r, t)$ with $\tilde u(r,t)$ in the region $(r, t)\in [h(t_0)-\delta^*, R_0]\times [t_0, t^0]$.
On the parabolic boundary of this region, by \eqref{*1}, \eqref{**1} and \eqref{***1}, we have $u\leq \hat u$, and hence we can use the comparison principle to conclude that
$u(r, t)\leq \hat u(r, t)$ in this region. In particular, $ u(\hat h(t^0), t^0)\leq \hat u(\hat h(t^0), t^0)=0$.
On the other hand, from $h(t^0)>\hat h(t^0)$ we obtain $u(\hat h(t^0), t^0)>0$. This contradiction proves our claim, and the proof of the lemma is now complete.
\end{proof}

\begin{Lemma}\label{subsequence}
There exists a sequence $\{t_k\}$ with $t_k\to\infty$ such that
\[
\liminf_{k\to\infty} u(h(t_k)-r, t_k)\geq \Phi_{c_*}(-r)>0 \mbox{ for every } r>0.
\]
\end{Lemma}
\begin{proof}
We divide the proof into two steps.

\noindent\textsc{Step 1. } We show that for each $r_0>0$, there exist $C_{r_0}<0$ and a sequence $t_n\to\infty$ such that
\[
\lim_{n\to\infty}u( h(t_n)-r, t_n)=\Phi_{c_*}(r_0+C_{r_0}-r) \mbox{ in } C_{\rm loc}(\mathbb R).
\]

 Let $\sigma(r)$ be defined as in Lemma \ref{sigma(r)}. For any given $r_0>0$, there exists a sequence
$t_n\to\infty$ satisfying
\[
\lim_{n\to\infty} u(h(t_n)-r_0, t_n)=\sigma(r_0)>0.
\]

From the definition of $\tilde u$ we obtain
\[
u(h(t)-r,t)=\tilde u(\tilde h(t)-r, t).
\]
In view of \eqref{C-3C} and Lemma \ref{lem:convergence.solution.subsequence}, by passing to a further subsequence of $\{t_n\}$ we may assume that $\tilde h(t_n)\to m\in [-C, C]$ and
\[
\lim_{n\to\infty}u( h(t_n)-r, t_n)=\lim_{n\to\infty}\tilde u(\tilde h(t_n)-r, t_n)=\Phi_{c_*}(C_0+m-r) \mbox{ in } C_{\rm loc}(\mathbb R)
\]
for some $C_0\in\mathbb R$.
Taking $r=r_0$ we obtain
\[
\Phi_{c_*}(C_0+m-r_0)=\sigma(r_0)>0 \mbox{ and hence } C_{r_0}:=C_0+m-r_0<0,
\]
and $\Phi_{c_*}(C_0+m-r)=\Phi_{c_*}(r_0+C_{r_0}-r)$. The proof of Step 1 is now complete.

\smallskip

\noindent \textsc{Step 2.} Find the required sequence $\{t_k\}$.

Let $r_k$ decrease to 0 as $k\to\infty$. By Step 1, for each $k\geq 1$, there exist $C_{r_k}<0$ and a sequence $\{t_n^{(k)}\}_{n=1}^\infty$
such that, as $n\to\infty$,
\[
t_n^{(k)}\to\infty,\; u(h(t_n^{(k)})-r, t_n^{(k)})\to \Phi_{c_*}(r_k+C_{r_k}-r) \mbox{ uniformly in } r\in [0, k].
\]
Since
\[
\Phi_{c_*}(r_k+C_{r_k}-r)>\Phi_{c_*}(r_k-r) \mbox{ for } r\in [r_k, k],
\]
 there exists $n_k$ large such that
\[
t^k_{n_k}>k,\; 
u(h(t_{n_k}^{(k)})-r, t_{n_k}^{(k)})>\Phi_{c_*}(r_k-r) \mbox{ for } r\in [r_k, k].
\]
Define $t_k:=t_{n_k}^{(k)}$. Then clearly
\[
\lim_{k\to\infty} t_k=\infty,\; \liminf_{k\to\infty}u(h(t_k)-r, t_k)\geq \lim_{k\to\infty}\Phi_{c_*}(r_k-r)=\Phi_{c_*}(-r)>0 \mbox{ for every } r>0.
\]
This completes the proof of the lemma.
\end{proof}

\begin{proof}[Proof of Proposition \ref{prop6.1}]
Let $\{t_k\}$ be the sequence given by Lemma~\ref{subsequence}. By Lemma~\ref{lem:convergence.solution.subsequence}, there is a subsequence of $\{t_k\}$, still denoted by itself, such that
\[
\lim_{k\to\infty}u(r+c_*t_{k}-(N-1)c^*\log t_{k},t_k)=\Phi_{c^*}(r-r_0)
\]
for some $r_0\in\mathbb{R}$. Moreover,  $\lim_{k\to\infty}\tilde{h}(t_k)=m\in [-C, C]$ and so necessarily $m\ge r_0$.

On the other hand, by Lemma \ref{subsequence} and our choice of $\{t_k\}$, as $k\to\infty$,
\[
u(h(t_k)-r, t_k)=\tilde u(\tilde h (t_k)-r, t_k)\to \Phi_{c_*}(m-r-r_0)\geq \Phi_{c_*}(-r).
\]
It follows that $m-r_0\leq 0$. Hence we must have $m=r_0$, and the proof is complete.
\end{proof}

\subsection{Proof of Theorem \ref{thm:main}.}

Let $\{t_n\}$ be a sequence having the properties stated in Proposition \ref{prop6.1}, and let $T, M_0, \bar M$ and $b$ be given in the proof of Lemma 3.4. For any given small $\epsilon>0$, we now choose $M_1>M_0$ large so that
\[
\frac {b\log M_1}{M_1}\in (0, \epsilon).
\]
Then define
\[
g_n(t):=1-\frac{\log(t-t_n+M_1)}{t-t_n+M_1},\; \underline h_n(t)=k(t)+\frac{b\log(t-t_n+M_1)}{t-t_n+M_1}+r_0-2\epsilon,
\]
and
\[
w_n(r,t):=g_n(t)\Phi\Big(r-\underline h_n(t);\alpha(c_*-\frac{(N-1)c^*}{t})\Big).
\]
We observe that for all large $n$, say $n\geq n_0$, we have
\[
\mbox{$T_n:=t_n+M_0-M_1>T$ and $t_0^n:=t_n-M_1\in (0, T_n-M_0]$.}
\]
Thus by the proof of Lemma 3.4 (with $t_0$ replaced by $t_0^n$ and $T$ replaced by $T_n$) we find that, whenever $n\geq n_0$,
\[
\mathcal L w_n\leq 0 \mbox{ for } t\geq T_n,\; r\in (c_*t-\bar M\log t, \underline h_n(t)),
\]
and
\[
w_n(r, t)\leq u(r,t) \mbox{ for } r=c_*t-\bar M\log t,\; t\geq T_n.
\]

We note that $T_n<t_n$ and so, if we can show, for large $n\geq n_0$,
\begin{equation}\label{t=t_k}
w_n(r, t_n)\leq u(r, t_n) \mbox{ for } r\geq c_*t_n-[\bar M-(N-1)c^*] \log t_n,
\end{equation}
then we can apply the comparison principle as before to conclude that
\begin{equation}\label{u-geq}
w_n(r, t)\leq u(r, t) \mbox{ for } t>t_n,\; r\geq c_* t-[\bar M-(N-1)c^*] \log t,
\end{equation}
which yields in particular
\[
h(t)\geq \underline h_n(t)=k(t)+\frac{b\log(t-t_n+M_1)}{t-t_n+M_1}+r_0-2\epsilon.
\]
It follows that
\[
\liminf_{t\to\infty}\; [h(t)-k(t)]\geq r_0-2\epsilon.
\]
Since $\epsilon>0$ can be arbitrarily small we thus obtain
\[
\liminf_{t\to\infty}\; [h(t)-k(t)]\geq r_0,
\]
provided  \eqref{t=t_k} holds.

We now prove \eqref{t=t_k}. For convenience, we write
\[
\tilde u_n(r):=u(r+k(t_n), t_n),\; \tilde w_n(r):=w_n(r+k(t_n),t_n).
\]
So \eqref{t=t_k} is equivalent to
\begin{equation}\label{t=t_k-1}
\tilde u_n(r)\geq \tilde w_n(r) \mbox{ for } r\geq -[\bar M-(N-1)c^*]\log t_n.
\end{equation}
From the definition of $w_n$ we have
\[
\tilde w_n(r)=\Big(1-\frac{\log M_1}{M_1}\Big)\Phi\left(r-\frac{b\log M_1}{M_1}-r_0+2\epsilon; \alpha \big(c_*-\frac{(N-1)c^*}{t_n}\big)\right).
\]
Hence, in view of $b(\log M_1)/M_1<\epsilon$,
\[
\lim_{n\to\infty} \tilde w_n(r)=\Big(1-\frac{\log M_1}{M_1}\Big)\Phi_{c_*}\left(r-\frac{b\log M_1}{M_1}-r_0+2\epsilon\right)
\leq \Big(1-\frac{\log M_1}{M_1}\Big)\Phi_{c_*}(r-r_0+\epsilon)
\]
uniformly for $r\in \mathbb R$.

On the other hand, by our choice of the sequence $\{t_n\}$, we have, for any $a\in\mathbb R$,
\[
\lim_{n\to\infty}\tilde u_n(r)=\Phi_{c_*}(r-r_0) \mbox{ uniformly in } r\in [a,\infty),\; \lim_{n\to\infty} [h(t_n)-k(t_n)]=r_0.
\]
Since $\Phi_{c_*}(-\infty)=1$, we may choose $a<0$ large negative so that
\[
\Phi_{c_*}(a-r_0)>\Big(1- \frac{\log M_1}{M_1}\Big)\Big/\Big(1- \frac{\log M_1}{2M_1}\Big).
\]
Then we can find $n_1\geq n_0$ large so that
\[
\tilde u_n(r)\geq \tilde w_n(r) \mbox{ for } r\geq a,\; n\geq n_1,
\]
where we have also used the facts that for all large $n$,
\[
\tilde u_n(r)\geq  \Big(1- \frac{\log M_1}{2M_1}\Big) \Phi_{c_*}(r-r_0)\geq \Big(1- \frac{\log M_1}{2M_1}\Big) \Phi_{c_*}(-\epsilon/2)>0 \mbox{ for } r\in [a, r_0-\epsilon/2]
\]
and $\tilde w_n(r)=0$ for $r\geq r_0-\epsilon/2$.

For $r\in \big[-[\bar M-(N-1)c^*]\log t_n, a\big]$, due to the monotonicity of $u(t,r)$ in $r\geq h(0)$, we have for large $n$,
\[
\tilde u_n(r)\geq \tilde u_n(a)\geq \Big(1- \frac{\log M_1}{2M_1}\Big) \Phi_{c_*}(a-r_0)\geq 1- \frac{\log M_1}{M_1},
\]
while
\[
\tilde w_n(r)\leq 1- \frac{\log M_1}{M_1}.
\]
Thus \eqref{t=t_k-1} holds, as we wanted. This proves \eqref{u-geq}.

From \eqref{u-geq} we obtain, for any $a\in\mathbb R$,
\[
\liminf_{t\to\infty} u(r+k(t), t)\geq \Phi_{c_*}(r-r_0+2\epsilon) \mbox{ uniformly for } r\geq a.
\]
Due to the monotonicity of $u(t,r)$ in $r\geq h(0)$, the facts $\Phi_{c_*}(-\infty)=1$ and $u(t,r)\to 1$ locally uniformly for
$r\in [0, \infty)$, and that $\epsilon>0$ can be arbitrarily small, it follows easily that
\[
\liminf_{t\to\infty} u(r+k(t), t)\geq \Phi_{c_*}(r-r_0) \mbox{ uniformly for } r\geq -k(t).
\]

By a parallel argument using the proof of Lemma 3.5, we can similarly show that
\[
\limsup_{t\to\infty}\; [h(t)-k(t)]\leq r_0\;\mbox{ and }
\]
\[
\limsup_{t\to\infty} u(r+k(t), t)\leq \Phi_{c_*}(r-r_0) \mbox{ uniformly for } r\geq -k(t).
\]
The conclusions in Theorem \ref{thm:main} thus all follow.\hfill $\Box$

\subsection{Remarks on the one dimension case and a gap in \cite{Biro-2002}} The nice paper of Bir\'o~\cite{Biro-2002} contains a gap in its proof of the main theorem in Section 4. It is proved in \cite{Biro-2002} that along a time sequence $t=\theta_{k}$, the free boundary $\xi=\zeta(t)$ of the solution $z(\xi, t)$ converges to some $\xi^*\in [\xi_1,\xi_2]$, and then along a subsequence $\theta_{k_s}$, $z(\xi,\theta_{k_s})\to w(\xi)$, and $w(\xi)$ is a shift
of the traveling wave $V(\xi)$, namely $w(\xi)=V(\xi-\hat \xi^*)$.  One expects that $\hat\xi^*=\xi^*$ and this fact
was used in \cite{Biro-2002} without a proof. This is the gap we would like to discuss here.

Let us now explain how this gap can be fixed.
Since $w(\xi)\equiv 0$ for $\xi\geq \xi^*$, necessarily $\hat\xi^*\leq \xi^*$. To show that $\hat\xi^*<\xi^*$ cannot happen, we need a certain non-degeneracy result for $z(\xi, t)$, which can be proved in the same way as in Lemmas~\ref{sigma(r)} and~\ref{subsequence} here, where the same non-degeneracy property is proved for our radially symmetric solution in $\mathbb R^N$. This approach  requires to choose the time sequence as in Section 6.2 here, differently from $t=\theta_k$ in \cite{Biro-2002}, and the subsolutions and supersolutions
 along the sequence $\{t_n\}$ are obtained from Lemma 3 of \cite{Biro-2002}, as in the proof of Lemma~\ref{G=const} here.  We leave the details to the interested reader.

With the above non-degeneracy property we could also extend the result of Theorem~\ref{thm:main} to the one-space dimension case with non-symmetric
initial function $u_0$, which is nonnegative and has non-empty compact support. The idea is to follow the approach here but replace the rather involved construction of subsolutions and supersolutions in Lemmas~3.4 and 3.5 here by corresponding ones in \cite{Biro-2002} (performed in $\mathbb R_+$ and in $\mathbb R_-$ respectively, instead of over $\mathbb R$, and combined with estimate at $x=0$ of the type given in Lemma 3.2 here), which are much simpler and do
not involve logarithmic corrections; see Section 3 there for details.  In this case the spatial support of $u(\cdot, t)$ for large times is an interval $[-h_-(t),h_+(t)]$, and the long-time dynamical behavior of the solution is given by the following theorem, with the detailed proof  left to the interested reader.

\begin{Theorem} Let $u_0(x)\not\equiv 0$ be nonnegative and compactly supported $($not necessarily symmetric$)$, and let  $u$ be the corresponding solution to~\eqref{eq:main} with $N=1$. Then there exist constants $r_\pm\in\mathbb{R}$, depending on $u_0$, such that $\lim_{t\to\infty}(h_\pm(t)-c_*t)=r_\pm$. Moreover,
$$
\begin{aligned}
&\lim_{t\to\infty}\left\{\sup_{x\in\mathbb{R}_+}\big|u(x,t)-U_{c_*}(x-c_*t-r_+)\big|\right\}=0,
\\
&\lim_{t\to\infty}\left\{\sup_{x\in\mathbb{R}_-}\big|u(x,t)-U_{c_*}(-x-c_*t-r_-)\big|\right\}=0.
\end{aligned}
$$
\end{Theorem}

\end{document}